\documentclass[reqno,12pt,a4paper]{amsart}
\UseRawInputEncoding
\usepackage{mathrsfs}
\usepackage{amssymb}
\usepackage{amsfonts}
\usepackage{latexsym}
\usepackage{amsthm}
\usepackage{graphicx}
\usepackage{wrapfig}
\def\lb{\label}
\newcommand{\er}[1]{\textrm{(\ref{#1})}}

\begin{document}
%%%%%%%%%% Some definitions %%%%%%%%%%
%%%%%%%% Equations, theorems %%%%%%%%%
\renewcommand{\theequation}{\arabic{section}.\arabic{equation}}
\theoremstyle{plain}
\newtheorem{theorem}{\bf Theorem}[section]
\newtheorem{lemma}[theorem]{\bf Lemma}
\newtheorem{corollary}[theorem]{\bf Corollary}
\newtheorem{proposition}[theorem]{\bf Proposition}
\newtheorem{definition}[theorem]{\bf Definition}
\newtheorem{remark}[theorem]{\bf Remark}
%%%%% Alphabet %%%%%
\def\a{\alpha}  \def\cA{{\mathcal A}}     \def\bA{{\bf A}}  \def\mA{{\mathscr A}}
\def\b{\beta}   \def\cB{{\mathcal B}}     \def\bB{{\bf B}}  \def\mB{{\mathscr B}}
\def\g{\gamma}  \def\cC{{\mathcal C}}     \def\bC{{\bf C}}  \def\mC{{\mathscr C}}
\def\G{\Gamma}  \def\cD{{\mathcal D}}     \def\bD{{\bf D}}  \def\mD{{\mathscr D}}
\def\d{\delta}  \def\cE{{\mathcal E}}     \def\bE{{\bf E}}  \def\mE{{\mathscr E}}
\def\D{\Delta}  \def\cF{{\mathcal F}}     \def\bF{{\bf F}}  \def\mF{{\mathscr F}}
\def\c{\chi}    \def\cG{{\mathcal G}}     \def\bG{{\bf G}}  \def\mG{{\mathscr G}}
\def\z{\zeta}   \def\cH{{\mathcal H}}     \def\bH{{\bf H}}  \def\mH{{\mathscr H}}
\def\e{\eta}    \def\cI{{\mathcal I}}     \def\bI{{\bf I}}  \def\mI{{\mathscr I}}
\def\p{\psi}    \def\cJ{{\mathcal J}}     \def\bJ{{\bf J}}  \def\mJ{{\mathscr J}}
\def\vT{\Theta} \def\cK{{\mathcal K}}     \def\bK{{\bf K}}  \def\mK{{\mathscr K}}
\def\k{\kappa}  \def\cL{{\mathcal L}}     \def\bL{{\bf L}}  \def\mL{{\mathscr L}}
\def\l{\lambda} \def\cM{{\mathcal M}}     \def\bM{{\bf M}}  \def\mM{{\mathscr M}}
\def\L{\Lambda} \def\cN{{\mathcal N}}     \def\bN{{\bf N}}  \def\mN{{\mathscr N}}
\def\m{\mu}     \def\cO{{\mathcal O}}     \def\bO{{\bf O}}  \def\mO{{\mathscr O}}
\def\n{\nu}     \def\cP{{\mathcal P}}     \def\bP{{\bf P}}  \def\mP{{\mathscr P}}
\def\r{\rho}    \def\cQ{{\mathcal Q}}     \def\bQ{{\bf Q}}  \def\mQ{{\mathscr Q}}
\def\s{\sigma}  \def\cR{{\mathcal R}}     \def\bR{{\bf R}}  \def\mR{{\mathscr R}}
                \def\cS{{\mathcal S}}     \def\bS{{\bf S}}  \def\mS{{\mathscr S}}
\def\t{\tau}    \def\cT{{\mathcal T}}     \def\bT{{\bf T}}  \def\mT{{\mathscr T}}
\def\f{\phi}    \def\cU{{\mathcal U}}     \def\bU{{\bf U}}  \def\mU{{\mathscr U}}
\def\F{\Phi}    \def\cV{{\mathcal V}}     \def\bV{{\bf V}}  \def\mV{{\mathscr V}}
\def\P{\Psi}    \def\cW{{\mathcal W}}     \def\bW{{\bf W}}  \def\mW{{\mathscr W}}
\def\o{\omega}  \def\cX{{\mathcal X}}     \def\bX{{\bf X}}  \def\mX{{\mathscr X}}
\def\x{\xi}     \def\cY{{\mathcal Y}}     \def\bY{{\bf Y}}  \def\mY{{\mathscr Y}}
\def\X{\Xi}     \def\cZ{{\mathcal Z}}     \def\bZ{{\bf Z}}  \def\mZ{{\mathscr Z}}
\def\O{\Omega}
\def\mb{{\mathscr b}}
\def\mh{{\mathscr h}}
\def\me{{\mathscr e}}
\def\mk{{\mathscr k}}
\def\mz{{\mathscr z}}
\def\mx{{\mathscr x}}
\newcommand{\gA}{\mathfrak{A}}          \newcommand{\ga}{\mathfrak{a}}
\newcommand{\gB}{\mathfrak{B}}          \newcommand{\gb}{\mathfrak{b}}
\newcommand{\gC}{\mathfrak{C}}          \newcommand{\gc}{\mathfrak{c}}
\newcommand{\gD}{\mathfrak{D}}          \newcommand{\gd}{\mathfrak{d}}
\newcommand{\gE}{\mathfrak{E}}
\newcommand{\gF}{\mathfrak{F}}           \newcommand{\gf}{\mathfrak{f}}
\newcommand{\gG}{\mathfrak{G}}           %\newcommand{\gg}{\mathfrak{g}}
\newcommand{\gH}{\mathfrak{H}}           \newcommand{\gh}{\mathfrak{h}}
\newcommand{\gI}{\mathfrak{I}}           \newcommand{\gi}{\mathfrak{i}}
\newcommand{\gJ}{\mathfrak{J}}           \newcommand{\gj}{\mathfrak{j}}
\newcommand{\gK}{\mathfrak{K}}            \newcommand{\gk}{\mathfrak{k}}
\newcommand{\gL}{\mathfrak{L}}            \newcommand{\gl}{\mathfrak{l}}
\newcommand{\gM}{\mathfrak{M}}            \newcommand{\gm}{\mathfrak{m}}
\newcommand{\gN}{\mathfrak{N}}            \newcommand{\gn}{\mathfrak{n}}
\newcommand{\gO}{\mathfrak{O}}
\newcommand{\gP}{\mathfrak{P}}             \newcommand{\gp}{\mathfrak{p}}
\newcommand{\gQ}{\mathfrak{Q}}             \newcommand{\gq}{\mathfrak{q}}
\newcommand{\gR}{\mathfrak{R}}             \newcommand{\gr}{\mathfrak{r}}
\newcommand{\gS}{\mathfrak{S}}              \newcommand{\gs}{\mathfrak{s}}
\newcommand{\gT}{\mathfrak{T}}             \newcommand{\gt}{\mathfrak{t}}
\newcommand{\gU}{\mathfrak{U}}             \newcommand{\gu}{\mathfrak{u}}
\newcommand{\gV}{\mathfrak{V}}             \newcommand{\gv}{\mathfrak{v}}
\newcommand{\gW}{\mathfrak{W}}             \newcommand{\gw}{\mathfrak{w}}
\newcommand{\gX}{\mathfrak{X}}               \newcommand{\gx}{\mathfrak{x}}
\newcommand{\gY}{\mathfrak{Y}}              \newcommand{\gy}{\mathfrak{y}}
\newcommand{\gZ}{\mathfrak{Z}}             \newcommand{\gz}{\mathfrak{z}}

\def\ve{\varepsilon} \def\vt{\vartheta} \def\vp{\varphi}  \def\vk{\varkappa}
\def\vr{\varrho}
\def\Z{{\mathbb Z}} \def\R{{\mathbb R}} \def\C{{\mathbb C}}  \def\K{{\mathbb K}}
\def\T{{\mathbb T}} \def\N{{\mathbb N}} \def\dD{{\mathbb D}} \def\S{{\mathbb S}}
\def\B{{\mathbb B}}
%%%%% Arrows %%%%%
\def\la{\leftarrow}              \def\ra{\rightarrow}     \def\Ra{\Rightarrow}
\def\ua{\uparrow}                \def\da{\downarrow}
\def\lra{\leftrightarrow}        \def\Lra{\Leftrightarrow}
\newcommand{\abs}[1]{\lvert#1\rvert}
\newcommand{\br}[1]{\left(#1\right)}
\def\lan{\langle} \def\ran{\rangle}
%%%%% Typography %%%%%
\def\lt{\biggl}                  \def\rt{\biggr}
\def\ol{\overline}               \def\wt{\widetilde}
\def\no{\noindent}
%%%%% Math signs %%%%%
\let\ge\geqslant                 \let\le\leqslant
\def\lan{\langle}                \def\ran{\rangle}
\def\/{\over}                    \def\iy{\infty}
\def\sm{\setminus}               \def\es{\emptyset}
\def\ss{\subset}                 \def\ts{\times}
\def\pa{\partial}                \def\os{\oplus}
\def\om{\ominus}                 \def\ev{\equiv}
\def\iint{\int\!\!\!\int}        \def\iintt{\mathop{\int\!\!\int\!\!\dots\!\!\int}\limits}
\def\el2{\ell^{\,2}}             \def\1{1\!\!1}
\def\sh{\sharp}
\def\wh{\widehat}
\def\bs{\backslash}
\def\na{\nabla}
\def\ti{\tilde}
%%%%% Math operations %%%%%
\def\sh{\mathop{\mathrm{sh}}\nolimits}
\def\all{\mathop{\mathrm{all}}\nolimits}
\def\Area{\mathop{\mathrm{Area}}\nolimits}
\def\arg{\mathop{\mathrm{arg}}\nolimits}
\def\const{\mathop{\mathrm{const}}\nolimits}
\def\det{\mathop{\mathrm{det}}\nolimits}
\def\diag{\mathop{\mathrm{diag}}\nolimits}
\def\diam{\mathop{\mathrm{diam}}\nolimits}
\def\dim{\mathop{\mathrm{dim}}\nolimits}
\def\dist{\mathop{\mathrm{dist}}\nolimits}
\def\Im{\mathop{\mathrm{Im}}\nolimits}
\def\Iso{\mathop{\mathrm{Iso}}\nolimits}
\def\Ker{\mathop{\mathrm{Ker}}\nolimits}
\def\Lip{\mathop{\mathrm{Lip}}\nolimits}
\def\rank{\mathop{\mathrm{rank}}\limits}
\def\Ran{\mathop{\mathrm{Ran}}\nolimits}
\def\Re{\mathop{\mathrm{Re}}\nolimits}
\def\Res{\mathop{\mathrm{Res}}\nolimits}
\def\res{\mathop{\mathrm{res}}\limits}
\def\sign{\mathop{\mathrm{sign}}\nolimits}
\def\span{\mathop{\mathrm{span}}\nolimits}
\def\supp{\mathop{\mathrm{supp}}\nolimits}
\def\Tr{\mathop{\mathrm{Tr}}\nolimits}
\def\BBox{\hspace{1mm}\vrule height6pt width5.5pt depth0pt \hspace{6pt}}
\def\where{\mathop{\mathrm{where}}\nolimits}
\def\as{\mathop{\mathrm{as}}\nolimits}
%%%%%%%%%%%%% specialities %%%%%%%%%%%%%%
\newcommand\nh[2]{\widehat{#1}\vphantom{#1}^{(#2)}}
%{{\mathop{#1}\limits^\wedge}\vphantom{#1}^{(#2)}}
\def\dia{\diamond}
\def\Oplus{\bigoplus\nolimits}
%%%%%%%%%%% End of definitions %%%%%%%%%%
%%%%% OLD OLD OLD
\def\qqq{\qquad}
\def\qq{\quad}
\let\ge\geqslant
\let\le\leqslant
\let\geq\geqslant
\let\leq\leqslant
\newcommand{\ca}{\begin{cases}}
\newcommand{\ac}{\end{cases}}
\newcommand{\ma}{\begin{pmatrix}}
\newcommand{\am}{\end{pmatrix}}
\renewcommand{\[}{\begin{equation}}
\renewcommand{\]}{\end{equation}}
\def\eq{\begin{equation}}
\def\qe{\end{equation}}
\def\[{\begin{equation}}
\def\bu{\bullet}
%\bigskip
%\medskip
%\smallskip

%\title[{Inverse resonance problem for Jacobi
%operators on half-lattice}] {Inverse resonance problem for Jacobi
%operators on half-lattice}

\title[{Inverse  problems  for Jacobi operators
% with finitely supported perturbations
}] {Inverse  problems  for Jacobi operators  with finitely supported
perturbations }

%Inverse  problems  for Jacobi operators  with finitely supported perturbations

\date{\today}
\author[Evgeny Korotyaev, Ekaterina Leonova]{Evgeny L. Korotyaev, Ekaterina Leonova}
\address{ Saint-Petersburg
State University, Universitetskaya nab. 7/9, St. Petersburg, 199034,
Russia and HSE University, 3A Kantemirovskaya ulitsa, St.
Petersburg, 194100, Russia, \ korotyaev@gmail.com, \
e.korotyaev@spbu.ru,\ eleonova@hse.ru }

\begin{abstract}
\no We solve the inverse  problem   for Jacobi operators on the half
lattice with finitely supported perturbations, in particular, in
terms of resonances. Our proof is based on the results for the
inverse eigenvalue problem for specific finite Jacobi matrices  and
theory of polynomials. We determine forbidden domains for resonances
and maximal possible multiplicities of real and complex resonances.
\end{abstract}
\maketitle

\section {Introduction and main results}
\setcounter{equation}{0}

\subsection{Introduction}
We consider Jacobi operators $J$ acting on $\ell^2(\N)$ and given
by
\[\lb{j}
(Jf)_x=a_xf_{x+1}+b_xf_x+a_{x-1}f_{x-1},\qqq x\in \N=\{1,2,3,...\},
\]
where $f=(f_x)_1^\iy\in \ell^2(\N)$ and formally $f_0=0.$ For two
real sequences $(a_x)_{1}^{\iy}$ and $b=(b_x)_{1}^{\iy}\in
\ell^\iy(\N)$ define a  perturbation of $J$ by
$$
q = q(a,b) = (b_1, a_{1}-1,
b_2, a_2-1,\dots).
$$
We assume that $q$ is finitely supported  and
belongs to the class $\gX_{k}, k\in \N$ given by
\[
\begin{aligned}
 \lb{qs}
&  \gX_k = \rt\{h\in \ell^{\infty}(\N): h_{2x}>-1
\;\forall x\in \mathbb{N} \text { and }  h_k\neq 0,\; h_x=0\;
\forall x>k \rt\}.
\end{aligned}
\]
Note that if $k\in\{2p-1, 2p\}$ for some $p\in \N$, then we have

if $x>p$, then $a_x-1=b_x = 0$,

if $k = 2p-1,$ then  $a_p=1, b_p\ne0$,

if $k = 2p,$ then $a_p\ne1$.

\no For a set $Y\subset\mathbb{Z}$ (below  $Y=\mathbb{N},\;
\mathbb{Z}_+=\N\cup \{0\}$ or $\mathbb{Z}$) define the spaces
$\ell^2(Y)$ and $\ell^{\infty}(Y)$ as the spaces of sequences
$(h_n)_{n\in Y}$ equipped with the norms
\[
\|h\|_2 = \Bigl(\sum_{n\in Y}|h_n|^2\Bigr)^{\frac{1}{2}}\geq 0,
\;\;\; \|h\|_{\infty} = \sup_{n\in Y}|h_n|
\]
respectively. Recall the main properties of $J,$ see, e.g.,
\cite{T89}. It has purely absolutely continuous spectrum $[-2,2]$
plus a finite number of simple eigenvalues on the set $\R\sm
[-2,2]$. Thus the spectrum of $J$ has the form
\[
\lb{sJ} \s(J)=\s_{ac}(J)\cup \s_d(J),\qqq \s_{ac}(J)=[-2,2],\qqq
\s_d(J)\ss \R\sm [-2,2],
\]
see Fig. 1. The eigenvalues of $J$ are denoted by
\[\l_{-1_-}< \dots< \l_{-n_-}<-2<2<\l_{n_+}<\dots< \l_{1},
\]
\setlength{\unitlength}{1.0mm}
\begin{figure}[h]
\centering
\unitlength 1.0mm % = 2.845pt
\begin{picture}(135,25)
\put(5,10){\line(10,0){100.00}}
\put(25,9){$\ts$} \put(25,13){$\l_{-1}$}
\put(51,9){$\ts$} \put(35,12.5){$\l_{-2}$}
\put(35,9){$\ts$}\put(51,12.5){$\l_{-3}$}
\put(56,6){-2} \put(58,10){\line(0,1){1.6}}
\put(58,11.6){\line(1,0){22.00}} \put(80,10){\line(0,1){1.6}}
\put(65,15){$\s_{ac}(J)$}
\put(78,6){2}
\put(85,9){$\ts$} \put(92,9){$\ts$} \put(85,12.5){$\l_{2}$}
\put(92,12.5){$\l_{1}$}
\end{picture}
%\vspace{-10mm}
\caption{\footnotesize The spectrum of $J$} \label{fig}
\end{figure}
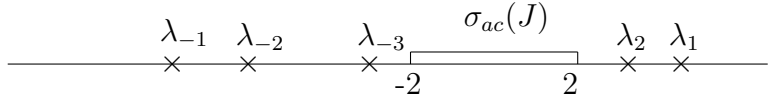
for some $n_\pm\ge 0$. The operator $J$ with a perturbation $q\in
\gX_{k},$ where $k\in\{2p-1, 2p\}$ has a representation as  a
semi-infinite matrix given by
\[
\label{J0} J=\left(\begin{array}{ccccccccc} b_1 & a_1 & 0 & 0 &  ...  &  0 & 0 & 0 & 0\cr
 a_1 & b_2 & a_2 & 0 &  ... &  0& 0 & 0 & 0 \cr
 0 & a_2 & b_3 & a_3 &  ... & 0 & 0 & 0 & 0\cr
 ... & ... & ... & ... & ... & ... & ... & ... & ... \cr
 0 & ...  & 0 & a_{p-1} & b_{p} & a_{p}& 0& 0 & 0\cr
 0 & ...  & 0 & 0 & a_{p} & 0 & 1 & 0 & 0\cr
 0 & ...  & 0 & 0 & 0 & 1 &0 & 1 & 0\cr
  ... & ... & ... & ... & ... & ... & ... & ... & ...
\end{array}\right).
\]
Define a new spectral variable $z = z(\l)$ by
\[\lb{lz}
\textstyle \l=\l(z)=z+{1\/z}, \qqq z\in \dD, \;\; \l\in \L:=\C\sm
[-2,2],
\]
where $\;\dD_r=\{z\in \C:|z|<r\}, r>0$ is the disc, and $\dD =
\dD_1$. Here $\l(z)$ is a conformal mapping from $\dD$ onto $\L$.
Its inverse function $z(\l)$ satisfies
\[
\textstyle{z(\l)={\l\/2}-\sqrt{{\l^2\/4}-1},\qqq
z(\l)=\frac{1+o(1)}{\l}\;\;\; \text{ as } \l\rightarrow\infty.}
\]
For the operator $J$ we define a finite Jacobi operator $J_p$ acting
on $\C^p$ and given by:
 \[
 \label{J1}
\Bigl(J_p f\Bigr)_x= a_xf_{x+1}+b_xf_x+a_{x-1}f_{x-1},\qqq
f=(f_n)_1^p\in \C^p,\qq
 x\in \N_p=\{1,2,\dots,p\},
\]
where formally $f_0=f_{p+1}=0.$ We define a matrix $J_p^1$ that is
obtained from $J_p$ by deleting the first row and the first column.
Similarly, the matrix $J_{p,1}$ is obtained from $J_p$ by deleting
the last row and the last column. If we delete both the first and
the last row and column, we obtain a matrix $J^1_{p,1}.$ The
operators corresponding to these matrices are self-adjoint.

The operator $J_p$ has $p$ simple eigenvalues labeled by
$\m_0<\m_1<..... <\m_{p-1}.$

The operator $J^1_{p,1}$ has $p-2$ simple eigenvalues  labeled by
$\n_1<\n_2<..... <\n_{p-2}.$

The operator $J_{p,1}$ has $p-1$ simple eigenvalues  labeled by
$\t_1<\t_2<..... <\t_{p-1}$.

The operator $J^1_p$ has $p-1$ simple eigenvalues labeled  by
$\vr_1<\vr_2<..... <\vr_{p-1}$.

\no The numbers $\t_j$ are  called Dirichlet eigenvalues, $\vr_j$
are called Neumann eigenvalues, $\m_j$ and $\n_j$ are called mixed
eigenvalues. Recall a well-known relation (see, e.g., \cite{vM76})
\[
\lb{bx}
\begin{aligned}
& \m_0<\ol {\varrho_1, \t_1}<\ol {\m_1,\n_1}<\ol {\varrho_2,
 \t_2}< \ol {\m_2,\n_2}<...<\m_{p-1},\\
\end{aligned}
\]
where $\ol {u,v}$ denotes $\min \{u,v\}\le \max \{u,v\}$.
 In order to describe eigenvalues of $J$ we use the set $n_\bu$ defined by
\[n_{\bu} = (-n_-,\dots, -1, 1, \dots, n_+).
\]
For a Jacobi operator $J$ we introduce the Jost solution
$(\p_x(z))_0^\iy$ of the equation
\[\lb{jej}
a_{x-1}\psi_{x-1}+a_x\psi_{x+1}+b_x\psi_x=(z+z^{-1}) \psi_x,\qqq \
|z|>0,\;\; x\in \Z_+,
\]
with initial conditions
 \[\p_x(z)=z^{x}\ \qq \ \forall  x>p.
\]
For the matrix $J$ we define the Jost function by $\p_0(z)$.  The
Jost function $\p_0(z)$ is a real polynomial of order $k.$ In
particular, we have
$$
\p_0(z) =\ca
    \frac{1}{A_0}+ O(z) & \as  \ z\to 0 \\
    Cz^{k}(1+o(1)) & \as \ \ z\to+\infty
  \ac,
$$
where the constants $A_0 = a_pa_{p-1}\cdot\dots\cdot a_1$ and $C\neq
0$ (see \cite{T89}).  Note that the function $\p_0$ has
$\gN=n_++n_-\ge 0$ simple zeros in $\dD$ given by
\[
\lb{eg1}
\begin{aligned}
\textstyle s_j = z(\l_j),\qq \l_j=s_j+{1\/s_j},\qqq j\in n_{\bu},
\\
 -1<s_{-n_-}<...<s_{-1}<0<s_{1}<...<s_{n_+}<1.
 \end{aligned}
\]
 The function $\p_0$ can also have a finite number of zeros in
$\C\sm \dD$ which are called resonances. The only possible zeros in
$\{|z|=1\}$ are $\pm 1,$ and these zeros are always simple. They are
called virtual states or resonances.

For two sequences $f=(f_x)_0^\iy, \ u =(u_x)_0^\iy$
we introduce the Wronskian by
\[\lb{wro}
\{f,u\}_x=a_x(f_x u_{x+1}-u_x f_{x+1}),\ \ x\ge 0.\] Note that if
$f$ and $u$ are some solutions of the Jacobi equation \er{jej}, then
$\{f,u\}_x$ does not depend on $x$. For each $J$ introduce finite
Jacobi matrices $J^{\pm}$ on $\mathbb{C}^p$ by
\[
\label{P1} J^{\pm} = J_p\pm a^2_p e_p^\top e_p,
\]
where the vector $e_p = (0,\dots,0, 1)\in \mathbb{C}^p.$ Each matrix
$J^{\pm}$ has $p$ simple eigenvalues $(\alpha^{\pm}_n)_1^p$. Below
we show that they satisfy the following relations, see Figure 2:
\[
\lb{ord}
\alpha^-_1< \alpha^+_1<\alpha^-_2<\alpha^+_2<...<\alpha^-_p<\alpha^+_p,
\]
\[
\lb{ord1}\alpha^-_1< \mu_0, \;\;\; \tau_{j} <\alpha^-_{j+1}
< \mu_j, \;\; j\in \mathbb{N}_{p-1},
\]
\[\lb{ord2}
\mu_{p-1}< \alpha^+_p,\;\;\; \mu_{j-1}<\alpha^+_j< \tau_{j}, \;\;
j\in \mathbb{N}_{p-1}.
\]

\subsection{Main results}
Let $\#(A, (c, d))$ denote the number of eigenvalues of an  operator
$A$ in the interval $(c, d)$, counting multiplicity. Now we describe
the location of the eigenvalues of the operators $J^{\pm}$ and $J,
J_p$ and $J_{p,1}$.

\begin{theorem}\lb{T1}
Let the operator $J$ be defined by (\ref{j}) with a perturbation $q
\in \gX_{k},$ where $k\in \{2p-1, 2p\}.$ Let
$\lambda_{-1}<\dots<\lambda_{-n_-}<-2$ and
$2<\lambda_{n_+}<\dots<\lambda_{1}$ be its eigenvalues. Then they
satisfy
\[
\lb{t13} n_- = \#(J^-, (-\infty, -2))\leq p,
\qq n_+ = \#(J^+, (2,+\infty))\leq p,
\]
\[\textstyle{\lb{t14}
n_-+n_+\leq\lfloor\frac{k}{2}\rfloor+1,}
\]
and
  \begin{multline}\label{t11}
    \alpha^-_1<\lambda_{-1}<\m_0<\alpha^+_1<\t_1<\alpha^-_2
    <\lambda_{-2}<\m_1<\dots \\
    \dots<\t_{n_- +1}<\alpha^-_{n_-}<\lambda_{-n_-}
    <\min\{-2, \alpha^+_{n_-}\} ,
  \end{multline}
  \begin{multline}\label{t12}
\max\{2, \alpha^-_{p - n_++1}\}
<\lambda_{n_+}<\alpha^+_{p-n_++1}<\t_{1}<\dots
\\
\dots<\m_{p-2}<\l_{2}<\alpha^+_{p-1}<\t_{p-1}<\alpha^-_p<\m_{p-1}
<\lambda_{1}<\alpha^+_p.
\end{multline}
Here $\lfloor x\rfloor$ is the floor function.
 In particular, if $k>1,$ then the operator $J$ has at least $p-1$
  resonances.
\end{theorem}
\begin{figure}
\begin{picture}(135,25)
\put(5,10){\line(1,0){120.00}}

\put(25,9){$\ts$} \put(25,13){$\a_1^-$}

 \put(51,9){$\ts$} \put(45,12.5){$\m_0$}
\put(45,9){$\ts$}\put(51,12.5){$\a_1^+$}

 \put(64,9){$\ts$} \put(71,9){$\ts$}
\put(64,13){$\t_1$} \put(71,13){$\a^-_2$}

\put(85,9){$\ts$} \put(92,9){$\ts$} \put(85,12.5){$\m_1$}
\put(92,12.5){$\a^+_2$}

\put(104,9){$\ts$} \put(110,9){$\ts$}
\put(104,13){$\t_2$} \put(110,13){$\a^-_3$}

\end{picture}\label{figu}
\caption{\footnotesize The eigenvalues of $J_p$ and $J^{\pm}$}
\end{figure}
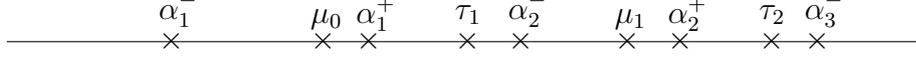
\textbf{Remark.} In the case of Schrodinger operators on the
half-line it is also possible to evaluate the number and location of
eigenvalues of the operator on the half-line using the eigenvalues
of the operator on the finite interval, see \cite{K21}. The
relations from \cite{K21} are similar to (\ref{t13}-\ref{t12}). But
in the continuous case the role of the eigenvalues
$(\a_j^{\pm})_1^p$ is played by the eigenvalues of the Sturm
-Liouville problem on a finite interval with mixed boundary
conditions, namely, Dirichlet condition on the one side and Neumann
condition on the other side. Thus, there is a principal difference
between the continuous and discrete cases.

\medskip

In Section 5 (Examples 1 and 2) we describe a location of
eigenvalues and resonances depending on the coefficients of the
matrix $J$ for $k = 1, 2$. The case $k = 1$ is trivial, while in the
case $k = 2$ there are several options for the arrangement of
eigenvalues and resonances. It is possible to give such description
for larger $k,$ but the number of different arrangements grows very
quickly.

For a real $t>0$ we define a set $\cE_p(t)$ of strictly increasing
sequences in $\mathbb{R}^p$ by
$$
\cE_p(t) = \left\{s = (s_k)_1^p\in\R^p: -t<s_1<s_2<\dots<s_p<t\right\}.
$$
Below we need a set $\cE_p = \cE_p(+\infty).$ For two interlacing
sequences $(x_n)_0^p$ and $(y_n)_1^p$ we define an increasing
sequence $x\star y$ by
\[
x\star y = (x_0, y_1,x_1,y_1,\dots,x_{p-1},y_p, x_p)\in \cE_{2p+1}.
\]
We introduce a set
$$
X_p = \gX_{2p}\bigcup\gX_{2p-1} = \rt\{h\in \ell^{\infty}(\N):
h_{2x}>-1 \;\forall x\in \mathbb{N},\;\;
|h_{2p-1}|+|h_{2p}|\neq 0,\; h_x=0\; \forall x>2p \rt\}.
$$
We prove some analogue of the Hochstadt result \cite{H78} (see
Theorem \ref{Ta}), which is used in the proof of main theorems.

\begin{proposition} \lb{TB}

Let the operator $J$ be defined by (\ref{j}) with a perturbation $q
\in X_p,\; p>0.$  Let $(\alpha^{\pm}_j)_1^{p}$ be the eigenvalues of
the operators $J^{\pm}.$ Then the mapping
$$\alpha^+\star\alpha^-: X_p\rightarrow \cE_{2p},$$
given by $q\mapsto \alpha^+\star\alpha^-$ is a a real analytic
isomorphism between $X_p$ and $\cE_{2p}$. Moreover, there is an
algorithm to recover $q$ from the numbers $\{\alpha^{\pm}_j\}.$
\end{proposition}

The resonances of $J$ can also be defined as the poles of  its
resolvent meromorphically continued in $\C$, see \cite{CK73}. The
S-matrix for $J$ is given by
\[
S(z)={\p_0(z^{-1})\/\p_0(z)}=e^{-i2\x(z)},\ \ \ \, \x(z)= \arg
\p_0(z), \ \ \ \ |z| = 1,
\]
where $\x$ is called the phase shift. The function $\p_0(z)$ has no
zeros in $\overline{\mathbb{D}}\setminus \Bigl([-1, s_{-1}]\bigcup
[s_1, 1] \Bigr).$ Therefore, we can choose a branch of $\log
\p_0(z)$ in this domain such that $\log \p_0(0)$ is real, which
uniquely defines the function $\x(z).$ The function $S$ has a
meromorphic continuation to the whole complex plane. Then the poles
of $S(z)$ are exactly the resonances of $J.$ Recall  that the
function $\x(z)$ is continuous in the circle $\{|z| =
1\}\setminus\{\pm1\}.$ It satisfies
\[\lb{lz1}
\textstyle
\x(1+0i)= -\pi n_++\frac{\pi m_+}{2}, \;\; \x(1-0i) =
 \pi \big (2n_- + n_+ + \frac{3m_+}{2}\big),
 \]
\[\lb{lz2}
\textstyle \x(-1\pm 0i) = \pi (n_-+n_+)\pm \frac{m_-}{2}\pi,
\]
where $m_{\pm}\leq 1$ is the multiplicity of $\pm 1$ as a zero of $\psi_0.$
 Define a set $X_p^+$ by
\[
X_p^+ = \gX^+_{2p-1}\bigcup\gX^+_{2p},
\;\text{ where }\;\gX_{k}^+=\rt\{q\in\gX_{k}:
\psi_0(z, q)\neq 0\;\;\forall z\in\mathbb{D}\rt\}.
\]
Recall that we have
$z(\l)={\l\/2}-\sqrt{{\l^2\/4}-1}, \;\;\l\in \Lambda.$
Below we consider a function $\x(z(\l+0i)), \; \l\in[-2,2].$
Note that for the conformal mapping $\l = z+\frac{1}{z}$ we have
$$
\begin{aligned}
z\Bigl([-2, 2]\pm 0i\Bigr) = \{z\in \C^{\mp}: |z|=1\}, \\
\x(z(\l+0i)) = - \x(z(\l-0i))\qqq \forall \qq \l\in[-2,2].
\end{aligned}
$$
We formulate the first results about inverse  problems  for Jacobi
operators  with finitely supported perturbations on the half
lattice.

\begin{theorem} \lb{pr1}

Let the operator $J$ be defined by (\ref{j}) with a perturbation $q
\in X_p^+$. Then the function $F(\l) = \frac{2\x(z(\l+0i))}{2p+1} +
\arg z(\l+0i)$ is strongly increasing in $\l\in (-2, 2),$ and each
equation $F(\l) = -\frac{\pi n}{2p+1},\; n\in\mathbb{N}_{2p} $ has a
unique solution $\omega_n \in (-2,2)$ which satisfies
\[
\omega_{2j-1} = \alpha^+_{p-j+1},\;\; \omega_{2j} = \alpha^-_{p-j+1}, \;\;\;
j\in \mathbb{N}_{p}.
\]
The mapping $\omega: X_p^+ \rightarrow
\cE_{2p}(2)$ given by
\[\lb{pr1w}
q\mapsto \omega = (\omega_n)_1^{2p}
\]
is a bijection between $X_p^+$ and $\cE_{2p}(2)$.
Moreover, there is an algorithm to recover the coefficients $a,b$
from the numbers $(\omega_j)_1^{2p}.$
\end{theorem}

\textbf{Remark.} 1) If $q\in X_p^+$ is given, then we compute $(\o_n)_1^{2p}\in \cE_{2p}(2)$, where $\o_{2j-1} = \a^+_{j}$ and $\o_{2j} = \a^-_{j}$. Conversely, if  $(\omega_n)_1^{2p}\in \cE_{2p}(2)$ is given, then we can compute the unique coefficients $a,b$.  In the proof of this theorem we only use well-known results on the inverse eigenvalue problem for specific finite Jacobi matrices and theory of polynomials. Marchenko equations are not used in the proof.

2) It is possible to use other eigenvalues (for example, $\m_j$ and $\t_j$) instead of $\alpha^{\pm}_j$ in this theorem. However, $\alpha^{\pm}_j$ are natural, since due to Theorem \ref{T1} the numbers $\alpha^{\pm}_j$ control the eigenvalues of $J$ and we just need  to stay that they belong to $[-2,2],$, which yields that $J$ has no eigenvalues. In the case with $\m_j, \;\t_j,$ we need more complicated conditions to be hold.

3) The recovering procedure is described in the proof of Proposition \ref{TB}. We also give an example of recovering the perturbation $q$ in the case $p = 2$ in Section 5, Example 6.
\medskip

Now we define the set characterising all possible  locations of the
eigenvalues and resonances.

\no {\bf  Definition R}. {\it Let $\cR_k\ss\C^k, k\in \N$ be a set of vectors $r=(r_n)_{1}^{k} $ such that:

\no R1) $0<|r_1|\le |r_2|\le \dots |r_k|,$

\no R2) a polynomial $f(z)=\prod_1^{k}(z-r_n) $ is real on the
real line,

\no R3) all zeros of $f$ in $\overline{\dD}$ are real and simple.
Denote the zeros of $f$ in $\mathbb{D}$ as $s_j,\; j \in n_{\bu}$ for
some $n_-, n_+\ge 0$ arranged by
$$
-1<s_{-n_-}<..<s_{-1}<0<s_{1}<...<s_{n_+}<1,
$$
and let
\[
 f(s^{-1}_j)\ne 0 \ \ \  \forall \ j \in n_{\bu},
\]
\no R4) $f$ has an odd number $\ge 1$ of zeros on each  interval
$(s_{j}^{-1}, s_{j-1}^{-1}),\; j = n_{\bu}\setminus\{-1\}$ and an
even number $\ge 0$ of zeros on each of intervals $(1,
s_{n_+}^{-1})$ and $(s_{-n_-}^{-1},-1)$.}
\begin{figure}
  \centering
  \includegraphics[scale=0.35]{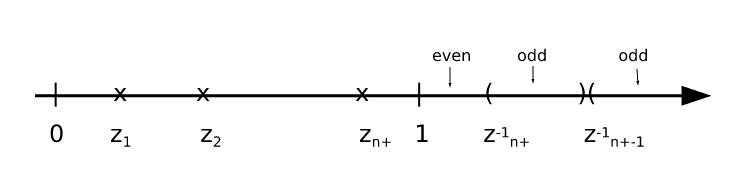}
  \caption{Position of the eigenvalues and resonances
   on the positive half-line}\label{R}
\end{figure}

The set of $\cR_k\ss\C^k, k\in \N$ is some analogue of the Jost
class functions for Schr\"odinger operator with a compactly
supported potential on the half-line from \cite{K04}.

 We sometimes write $\psi_n(q), r_{n} (q),..$ instead of
$\psi_n,r_{n}, ..$ when several potentials are being dealt with. Now
we construct the mapping $\gr: \gX_k\to \cR_k,\; k\in\mathbb{N}$ by
\[
 q \to \gr(q)=(r_n(q))_1^k\in \cR_k,
 \]
 where $r=(r_n)_1^k$ is the sequence of roots of the corresponding
  Jost function $\psi_0(z).$
\begin{theorem}\lb{T4}
Each mapping $\gr :  \gX_k\to \cR_k,\; k\in \mathbb{N}$ is a bijection between $\gX_k$ and $\cR_k.$  There is an
algorithm to recover the perturbation $q$ from
$\gr(q)$.\end{theorem}

\textbf{Remark.} This theorem solves the inverse resonance problem
for Jacobi operators on the half-lattice with  finitely-supported
perturbations. In our proof we do not use a  Marchenko equations
from  \cite{C73} to show  that the mapping $\gr$ is a surjection,
like it is usually done. Our proof consists of two parts. Firstly we
consider the case when the Jost function of $J$ does not have any
zeros in unit circle and prove Theorem \ref{pr1}. Secondly we prove
by induction that we can add roots to the unit circle. We also
present the algorithm to recover the matrix $J$ from its eigenvalues
and resonances. We only use theory of polynomials and some simple
lemmas proved beforehand in the second part.

\medskip

We also prove a corollary of Theorem \ref{T4} which shows that the
operator $J$ is uniquely determined by the values of the phase shift
function $\xi$ at a finite number of points.

\begin{corollary}\lb{corxi}

Consider two Jacobi operators $J_1$ and $J_2$ with perturbations
$q_1, q_2\in \gX_k,\; k\geq 1.$ Let $\xi_1$ and $\xi_2$ be the
corresponding phase shift functions. Let $(w_j)_1^{k}$ be a sequence
of distinct numbers on $\mathbb{C}_-\bigcap \{|z| = 1\}$ such that
$\xi_1(w_j) = \xi_2(w_j)$ for all $j\in\mathbb{N}_{k}.$ Then we have
$J_1 = J_2.$\end{corollary}

\textbf{Remark.} 1) One can see from the proof that this corollary
stays true if the numbers $(w_j)_1^{k}$ are chosen from the set
$\{|z| = 1\}\setminus \{\pm1\}$ such that $w_i \neq  \overline{w_j}$
for all $i,j\in \mathbb{N}_k.$

2) For the uniqueness we use a finite sequence $\{\x(w_j)\}_1^k $
for some $(w_j)_1^k$. Note that the number of parameters is equal to
the number of perturbations in $J.$ We do not need eigenvalues and
norming constants, like in general theory, see, e.g., \cite{CK73}.

\medskip

Below we discuss the location of resonances and specify the
forbidden domain.
\begin{theorem}\lb{cr1}

Let $\psi_0$ be the Jost function for some operator $J$ defined by
(\ref{j}) with a perturbation $q\in \gX_{k}.$ Then all  real
resonances of multiplicity greater than 1 and non-real resonances
belong to the disc $\{|z|<R_o\},$ where $R_o$ is given by:

if $k=2p \; (i.e.,\;a_p\neq 1),$ then

\[ \lb{cor61}
\textstyle
 R^2_o = \frac{1}{\beta{|1 - a_p^2|}},\;\;
\b = \ca 1, & {\rm if}\ \   a_p<1 \\
 |s_1s_{-1}|, & {\rm if}\ \  a_p>1
 \ac,
\]
if $k=2p-1\; (i.e.,\;a_p = 1, b_p\neq0),$ then
\[ \lb{cor62}
\textstyle R^2_o = \frac{1}{\beta|b_p|}, \;\;\b = \ca
|s_{-1}|, & {\rm if }\ \  b_p<0 \\
s_1, & {\rm if }\ \  b_p>0 \ac.
\]
Here $s_1>0$ and $s_{-1}<0$ are the roots of $\psi_0$ in
$\mathbb{D}$  with the smallest absolute value, see \er{eg1}. If
there are none of those, they are replaced with 1.

Moreover, these estimates are sharp: if $a_p$ (or $b_p$) and
$s_{\pm1}$ are fixed, then one of the resonances can be arbitrarily
close to the boundary of the disc for some Jacobi operator $J.$
\end{theorem}
Now we determine the maximal value of positive simple resonances.
\begin{theorem}
Let $\psi_0$ be the Jost function for some operator $J$ defined  by
(\ref{j}) with a perturbation $q\in \gX_{k}.$ Then all its positive
resonances of multiplicity 1 belong to the interval $(1, R_+),$
where $R_+$ is defined by:

if $k=2p \; (i.e.,\; a_p\neq 1),$ then

\[
\textstyle
 R_+ = \frac{1}{\beta_+|1 - a_p^2|},\;\; \;\;
\b_+ = \ca s_1, & \ \ {\rm if}\ \  a_p<1 \\
         \max\{\frac{s_1}{a_p^2-1}, |s_{-1}|\}, & \ \ {\rm if}\ \  a_p>1
\ac,
\]
if $k=2p-1\; ( i.e.,\; a_p = 1, b_p\neq0),$ then
\[
R_+ = \frac{1}{\b_+|b_p|}, \;\;
\b_+ = \ca |s_1s_{-1}|, & \ \ {\rm if}\ \  b_p<0 \\
           \max\{\frac{s_1}{|b_p|}, 1\}, & \ \ {\rm if}\ \ b_p>0\ac.
\]
Here $s_1>0$ and $s_{-1}<0$ are the roots of $\psi_0$  in
$\mathbb{D}$ with the smallest absolute value, see \er{eg1}. If
there are none of those, they are replaced with 1.

Moreover, these estimates are sharp: if $a_p$ (or $b_p$) and
$s_{\pm1}$ are fixed, then one of the resonances can be arbitrarily
close to $R_+$ for some Jacobi operator $J.$
\end{theorem}
\textbf{Remark.} 1) Similar result is proved for negative
resonances in Corollary \ref{corneg}.

2) For Schrodinger operators on the half line, the forbidden domain
for resonances is determined by some logarithmic curve, see e.g.,
\cite{K04}.

3) Thus, we specified the forbidden domain for the resonances  of
the operator $J.$ Namely, we obtained that if $k$ and $a_p$ (or
$b_p$) and the smallest eigenvalue are fixed, that all non-real
resonances and resonances of multiplicity greater than 1 are located
inside a ring of a fixed radius, while simple real resonances belong
to two fixed intervals, see Fig. 4.

\begin{figure}
\includegraphics[width=0.65 \textwidth]{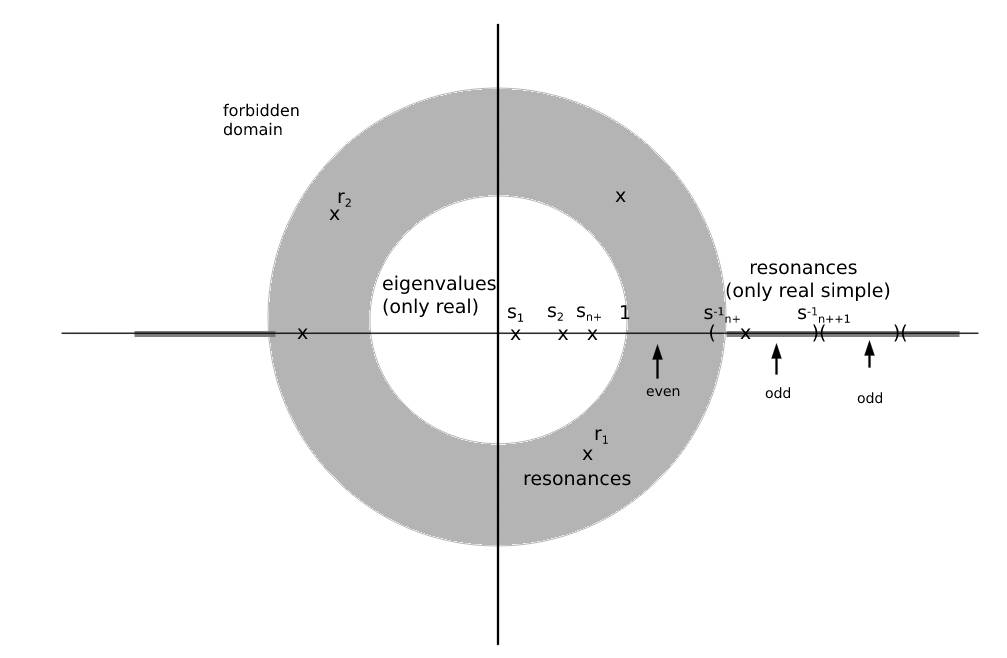}
\caption{the location of resonances on the complex plane}\label{est}
\end{figure}

\medskip

Now we discuss the maximal possible multiplicity of real-valued and
complex resonances.

\begin{corollary}\lb{crl}

Let $J$ be a Jacobi operator defined by (\ref{j}) with a
perturbation $q\in \gX_{k}.$ Let its Jost function $\psi_0$ have
$\gN = n_+ + n_-$ roots in $\mathbb{D}$ and let $r_o$ be a resonance
of $J.$ Then the maximal possible multiplicity of $r_o$ equals $M,$
where $M$ is given by:

i) if $k = 2p$ and $a_p<1$, then
\[
M =  \left\{
     \begin{array}{ll}
       p - \gN, & \hbox{$\text{if } r_o\in \mathbb{C}\setminus\mathbb{R}$} \\
       \min(2p - 2\gN + 1, 2p), & \hbox{if $r_o\in \mathbb{R}$}
     \end{array}
   \right.
,
\]

ii) if $k = 2p$ and $a_p>1$, then
\[ M = \left\{
     \begin{array}{ll}
       p+1 - \max(n_+, 1) - \max(n_-, 1), & \hbox{$\text{if } r_o\in \mathbb{C}\setminus\mathbb{R}$} \\
       2p + 3 - 2\max(n_+,1) -2\max(n_-,1), & \hbox{if $r_o\in \mathbb{R}$}
     \end{array}
   \right. ,
\]

iii) if $k = 2p-1$ and $b_p>0$, then
\[M =  \left\{
     \begin{array}{ll}
       p - \max(n_+, 1) - n_-, & \hbox{$\text{if } r_o\in \mathbb{C}\setminus\mathbb{R}$} \\
       2p + 1 - 2\max(n_+,1)-2n_-, & \hbox{if $r_o\in \mathbb{R}$}
     \end{array}
   \right. ,
\]

iv) if $k = 2p-1$ and $b_p<0$, then
\[ M = \left\{
     \begin{array}{ll}
       p - \max(n_-, 1) - n_+, & \hbox{$\text{if } r_o\in \mathbb{C}\setminus\mathbb{R}$} \\
       2p + 1 - 2\max(n_+,1)-2n_+, & \hbox{if $r_o\in \mathbb{R}$}
     \end{array}
   \right. .
\]
A complex resonance of the maximal multiplicity can  be located in
any  point in
$\mathbb{C}\setminus(\mathbb{R}\bigcup\overline{\mathbb{D}})$. If
$s_j, \; j\in n_{\bu}$ are the roots of $\psi_0$ in $\mathbb{D},$
then a real resonance of the maximal multiplicity can be located in
any of the intervals $(s_j^{-1}, s_{j-1}^{-1}), \; j=-n_-+1,\dots,
-1, 2,\dots, n_+.$
\end{corollary}

\textbf{Remark.} Thus roughly speaking the theorem gives that the
maximal multiplicity of resonaces $M=p-n_--n_+$. It follows from
Theorems \ref{T1} and  \ref{T4} that the numbers $n_+$ and $n_-$ can
be arbitrary, under the condition $n_- + n_+ \leq \lfloor
\frac{k}{2} \rfloor+1.$

\medskip

Finally we describe the behavior of the coefficients of the operator
$J$ when one of resonances moves to infinity.

\begin{corollary}\lb{cormax}
Let the operator $J$ with a perturbation $q\in \mathfrak{X}_{k}$ be  defined by
(\ref{j}). Let $\{r_n\}_{1}^{k}$ be the set of roots of its Jost function and let $(s_j), \; j\in n_{\bu}$ be the set of its eigenvalues.
We introduce the pair $(r, \gs)$ which satisfies one of the following conditions:

i) $r$ is some of non-real resonances and $\gs = \bar{r},$

ii) $r$ is some resonance from the interval $(s_{j+1}^{-1}, s_{j-1}^{-1})$ for some $j\in\{2, 3,\dots, n_+-1\},$ or $r$ is some resonance from the interval $(1, s_{n_+-1}^{-1})$ and $j = n_+.$ In both cases $\gs = s_j.$

%i) Let $\{r_n\}_{1}^{k}$ be the set of roots of its Jost
%function and let $r^o$ be some of its non-real resonances. Then the  resonance $r^o$ can be moved to any point in  $\wt r\in \C\sm
%(\R\bigcup\overline{\mathbb{D}})$ (preserving the conjugacy), and
%there exists a matrix $\wt J$ with $\wt q\in\mathfrak{X}_{k},$
%corresponding to the obtained set of roots. There also exists a
%matrix $J_o$ with a perturbation $q_0(a^o,
%b^o)\in\mathfrak{X}_{k-2}$ which Jost function has the roots
%$\{r_n\}_1^{k}\setminus\{r^o,\overline{r^o}\}.$

%In particular, if $k = 2p\; (i.e.,\;a_p\neq1),$ then we have
%\[
%a_p = 1 - \frac{1-(a^o_{p-1})^2+o(1)}{2|\wt r|^2} \;\;\;
%\text {as } |\wt r|\rightarrow +\infty.
%\]
%If $k = 2p-1\; (i.e.,\;a_p = 1, b_p \neq 0),$ then we have
%\[b_p = \frac{b^o_{p-1}}{|\wt r|^2} \rightarrow 0  \;\;\; \text
% {as } |\wt r|\rightarrow +\infty.\]
%Moreover, in both cases $\wt q\rightarrow q_o$  as $|\wt
%r|\rightarrow\infty.$

%ii) Let the operator $J$ have $n_+$ positive eigenvalues $(s_j)_1^{n_+}.$ Let $1<j< n_+$ and let $r$ be some resonance lying in the interval $(s_{j+1}^{-1}, s_{j-1}^{-1})$ or let $j = n_+$ and let $r$ be some resonance from the interval $(1, s_{n_+-1}^{-1})$.

Then the pair $(r, \gs)$ can be moved to any pair $(r_{\bu}, \gs_{\bu})$ on $\R\setminus
[-s_{-1}^{-1}, s_1^{-1}]$ or any two conjugate points in $ \C\sm
(\R\bigcup\overline{\mathbb{D}})$, and
there exists a matrix $\wt J$ with $\wt q\in\mathfrak{X}_{k},$
corresponding to the obtained set of roots. There also exists a
matrix $J_o$ with a perturbation $q_0(a^o,
b^o)\in\mathfrak{X}_{k-2}$ which Jost function has the roots
$\{r_n\}_1^{k}\setminus \{r, \gs\}.$

In particular, if $k = 2p\; (i.e.,\;a_p\neq1),$ then we have
\[
a_p = 1 - \frac{1-(a^o_{p-1})^2+o(1)}{2|\ga|} \;\;\;
\text {as } |\ga|\rightarrow +\infty,
\]
where $\ga =r_{\bu} \gs_{\bu} .$ If $k = 2p-1\; (i.e.,\;a_p = 1, b_p \neq 0),$ then we have

\[b_p = \frac{b^o_{p-1}}{|\ga|} \rightarrow 0  \;\;\; \text
 {as } |\ga|\rightarrow +\infty.\]
Moreover, in both cases $\wt q\rightarrow q_o$  as $|\ga|\rightarrow\infty.$
\end{corollary}
\textbf{Remark.} 1) While in the first case the points $r$ and $\gs$ can be moved to infinity only together, preserving the conjugacy, in the second case we can move both points or fix one of them and move another one to infinity, and still have $|\ga|\rightarrow +\infty.$

2) Similar result holds true for negative eigenvalues. We discuss analogous results for real resonances and the eigenvalues of the smallest absolute value in Corollary \ref{cormax1}.

3) Note that we cannot move the eigenvalues out of the circle without changing the location of resonances, as it will break the condition R4.

4) Similar problems were considered by Marletta and Weikard in \cite{MW07}. They only studied the case when all roots of the Jost function have a large absolute value. Our results describe the change of the matrix $J$ when only one or two roots are moved to infinity. After applying Corollary \ref{cormax} enough times, we obtain the result from \cite{MW07}.

%\medskip

\subsection{Short review}

A lot of papers are devoted to resonances of  one-dimensional
Schr\"odinger operators with compactly supported potentials, see
Froese \cite{F97}, Hitrik \cite{H99}, Korotyaev \cite{K04}, Simon
\cite{S00}, Zworski \cite{Z87} and references therein.
 Inverse problems (uniqueness, reconstruction, characterization) in
 terms of resonances were solved by Korotyaev for a
Schr\"odinger operator with a compactly supported potential on the
real line \cite{K05} and the half-line \cite{K04}. See also Zworski
\cite{Z01}, Brown-Knowles-Weikard \cite{BKW03} concerning the
uniqueness. The resonances for periodic plus a compactly supported
potential were considered  by Firsova \cite{F84}, Korotyaev
\cite{K11a}, Korotyaev-Schmidt \cite{KS12}. Christiansen \cite{C06}
discussed resonances for steplike potentials. The "local resonance"
stability problems were considered in \cite{K04a}, \cite{K16},
\cite{MSW10}. Resonances for three and fourth order differential
operator with compactly supported coefficients on the line were
studied in \cite{K19}, \cite{BK19}. Resonances for Stark operator
with compactly supported potentials were discussed in \cite{FH21},
\cite{K17}. Inverse resonance problems for Dirac operators with a
compactly supported potential on the half line is solved in
\cite{KM21}.

We briefly describe the results on the resonances for finitely
supported perturbations of 1dim discrete Laplacian. Inverse problem
on the real line was solved by Korotyaev \cite{K11a}, see also
Brown-Naboko-Weikard \cite{BNW05} about the uniqueness. Bledsoe
\cite{B12}, Marletta and Weikard \cite{MW07} studied stability for
Jacobi operators. Uniqueness and stability are also discussed in
\cite{MNSW12}. The case of periodic Jacobi operators with finitely
supported perturbations was studied by Iantchenko and Korotyaev
\cite{IK11}, \cite{IK12A}, \cite{IK12B} and Kozhan \cite{Ko16}. The
scattering problem for Jacobi operators with matrix-valued
coefficients was studied by Aptekarev and Nikishin \cite{AN83}.

\medskip

The plan of this paper is as follows. In Sect. 2 we discuss
the standard facts about Jost functions and fundamental solutions.
In Sect. 3 we prove the main results about the inverse problem. In
Sect. 4 we prove the main results about a location of resonances. In
Sect. 5 we consider some examples.

\section {Preliminaries}
\setcounter{equation}{0}

In order to prove the main theorems, we need the following asymptotics
\begin{lemma}\lb{L1}
\lb{f01} Let $q(a,b)\in \gX_k,\; k\in\{2p,2p-1\}$. Then

\no i) Each function $\p_{p-n}(z), n=0,...,p$ is a real polynomial
and satisfies
\[
\lb{f01-1} \p_{p}(z)={z^{p}\/a_p},\qqq
 \p_{p-n}(z)={z^{p+n}\/A_{p-n}}\rt(c_p-{c_pB_{n-1}+b_p\/z}+O(z^{-2})\rt),
 \qq
\]
\[
\lb{f01-2}
 \p_{0}(z)={z^{2p}\/A_0}\rt(c_p-{c_pB_{p-1}+b_p\/z}+O(z^{-2})\rt),
\]
as $z\to \iy$, where $B_j=\sum_{s=1}^{j}b_{p-s},
A_x=a_pa_{p-1}\cdot \cdot a_{x}$ and $c_x = 1 - a_x^2.$ Moreover,
\[
\lb{f01-3} \p_{p-n}(z)={z^{p-n}\/A_{p-n}}\rt(1-z\sum_{p-n+1}^p
b_k+O(z^{2})\rt),\qqq \p_{0}(z)={1+B_p z+O(z^2)\/A_0}\qq as \ z\to0.
\]
\no ii) The Jost function $\p_0(z)$ has the form:
\[
\lb{f02-1}
 \p_0(z)={1-a_p^2\/A_0}\prod_1^{2p}(z-r_n),\qqq
 \prod_1^{2p}r_n={1\/1-a_p^2},\qqq
 if \ a_p\ne 1,
\]
\[
\lb{f02-2} \p_0(z)=-{b_p\/A_0}\prod_1^{2p-1}(z-r_n),\qqq
 \prod_1^{2p-1}r_n={1\/b_p},\qq if \ a_p=1, b_p \neq 0,
 \qq
\]
\[
\lb{f02-3} \sum_1^{k}\frac{1}{r_n}=\sum_1^{p}b_n.
 \qq
\]
\end{lemma}
{\no\bf Proof.} i) Using the equation
$a_{x-1}f_{x-1}=(\l-b_x)f_x-a_xf_{x+1}$ and $\p_x=z^x,\; x>p$, we obtain:
\no if $n=0$, then
\[
\lb{asp0} \p_p={(\l-b_{p+1})\p_{p+1}-\p_{p+2}\/a_p}
={(z+z^{-1})z^{p+1}-z^{p+2}\/a_p}={z^{p}\/a_p},
\]
if $n=1$, then
\[
\lb{asp1} \p_{p-1}={(\l-b_{p})\p_{p}-a_p\p_{p+1}\/a_{p-1}}
={(z+z^{-1}-b_{p})z^{p}-a_p^2z^{p}\/a_pa_{p-1}}
={z^{p+1}c_p-b_pz^{p}+z^{p-1}\/a_pa_{p-1}},
\]
if $n=2$, then
\begin{multline}
\lb{asp2} \p_{p-2}={(\l-b_{p-1})\p_{p-1}-a_{p-1}\p_{p}\/a_{p-2}}
={(z+z^{-1}-b_{p-1})(z^{p+1}c_p-b_pz^{p}+
z^{p-1})-a_{p-1}^2z^{p}\/a_pa_{p-1}a_{p-2}}
\\
=z^p{z^{2}c_p-z(b_p+c_pb_{p-1})+(c_p+c_{p-1}+b_pb_{p-1})
-z^{-1}(b_p+b_{p-1})+z^{-2}\/a_pa_{p-1}a_{p-2}}
\\=
{z^{p+2}\/a_pa_{p-1}a_{p-2}}\rt(c_p-{b_p+
c_pb_{p-1}\/z}+{c_p+c_{p-1}+b_pb_{p-1}\/z^2}
-{b_p+b_{p-1}\/z^3}+{1\/z^4}\rt).
\end{multline}
Repeating this procedure we obtain \er{f01-1}, \er{f01-2}.

ii) Asymptotic \er{f01-1} yield \er{f02-1}. The identities
\er{f02-1} and $\p_{0}(0)={1\/A_0}>0$  (see \er{f01-3}) give
\er{f02-2}. \er{f02-3} follows from \er{f01-3} and \er{f02-2}. \BBox

For each eigenvalue $\l_j$ of $J$ we define a norming constant $\gc_j$ by
\[
\gc_j=\sum_{x\ge 1} \p_x(s_j)^2>0.
\]
 There exists a following identity, see e.g., [Ca]:
\[\lb{wl1}
\gc_j=a_1\p_0'(s_j)\p_1(s_j)z'(\l_j), \ \ \ \ \ s_j\in (-1,1).
\]
Let $\# (f, I)$ denote the number of zeros of $f$ on the set $I$. We
describe the location of eigenvalues and  real resonances of the
operator $J.$ We will use it to prove Theorems \ref{T1} and
1.5--1.6.
\begin{lemma}
\lb{f02} Let $q\in \gX_k,\; k\in\{2p,2p-1\}$ and let $-1<s_{-n_-}
<...<s_{-1}<0<s_{1}<...<s_{n_+}<1$ be all
zeros of $\p_{0}(z)$ on $(-1,1).$ Let $s^{-1}_j = \wh s_j.$
Then we have
\[
\lb{f01-4}
 \gc_j=a_1z'(\l_j)\p_0'(s_j)\p_1(s_j)
 =s_ja_1{\p_0'(s_j)\/ \p_0(\wh s_j)}>0, \ \ \forall j\in n_{\bu},
\]
\[
\lb{f01-5} (-1)^{j}\p_0'(s_j)>0,\ \ \
s_j(-1)^{j}\p_0(\wh s_j)>0, \ \ \forall j\in n_{\bu},
\]
\[
\lb{3.7}  \# (\p_0, (\wh s_{j+1},\wh s_{j}))=
{\rm odd}\ge 1,\;\;\; \forall j\in\mathbb{N}_{n_+-1},\ \ \
\]
\[
\lb{3.8} \# (\p_0, (\wh s_{j},\wh s_{j-1}))=
{\rm odd}\ge 1,\;\;\; \forall j\in [-n_- + 1: -1],
\]
\[
\lb{3.81} \# (\p_0, (s_{n_+},\wh s_{n_+}))={\rm even}\ge  0, \;\;\;
\#(\p_0, (\wh s_{-n_-}, s_{n_-}))={\rm even}.
\]
Moreover, if $ k=2p$ $(i.e.,\; a_p\neq 1)$, then we have
\[\lb{cpinf}
\ca a_{p}<1\qqq \Rightarrow  \qqq \# (\wh s_{\pm 1},\pm \iy)
={\rm odd}\ge 1,\qq \# (\wh s_{\pm n_\pm},\pm \iy)\ge n_\pm
\\
a_{p}>1 \qqq  \Rightarrow \qqq  \# (\wh s_{\pm 1},\pm \iy)
={\rm even}\ge 0,\qq \# (\wh s_{\pm n_\pm},\pm \iy)\ge n_\pm-1 \ac,
\]
if \ $k=2p-1$ $(i.e.,\; a_p = 1, b_p\neq 0)$, then we have
\[
\ca b_{p}>0\qqq \Rightarrow  \qqq \# (\wh s_1,\iy)=
{\rm odd}\ge 1,\qq \# (\wh s_{n_+},+ \iy)\ge n_\pm
\\
b_{p}<0 \qqq  \Rightarrow \qqq  \# (\wh s_1,\iy)=
{\rm even}\ge 0 ,\qq \# (\wh s_{n_+},+ \iy)\ge n_+-1 \ac.
\]
\end{lemma}
\no{\bf Proof.} Denote $f^{\diamond}(z) = f(\frac{1}{z}).$ Computing
the  Wronskian for the functions $\p$ and $\p^\diamond,$ we obtain
$$\{\p, \p^\diamond\}_0=\p_0^{\diamond} \p_1 - \p_1^{\diamond} \p_0,$$
$$\{\p, \p^\diamond\}_0 = \{\p, \p^\diamond\}_{p+1} =
 \p_{p+1}^{\diamond} \p_{p+2} -
 \p_{p+2}^{\diamond} \p_{p+1} = z - \frac{1}{z} = z\l'(z).$$
Then at $z = s_j$ we have
\[\lb{wl2}
\p_1(s_j) = \frac{s_j \l'(s_j)}{\p_0^{\diamond}(s_j) } ,
\]
and, substituting this identity into \er{wl1}, we obtain \er{f01-4}.

The function $\p_0(z),\; z\in (-1,1)$ has the zeros  $\{s_j\},$ and
$\p_0(0)>0$. Then we have $\p_0'(s_1)<0$ and $\p_0(s_2)>0$ etc.
Hence due to \er{f01-4} we obtain \er{f01-5}. Other identities
follow directly from this. \BBox

Now let us study the relationship between the eigenvalues of  the
operator $J$ and those of the finite Jacobi matrix $J_p$. Introduce the
fundamental solutions $\varphi = (\varphi_x(\l))_0^{\infty}$ and
$\vartheta = (\vartheta_x(\l))_0^{\infty}$ of the equation
\[\lb{je}
a_{x-1}f_{x-1}+a_xf_{x+1}+b_xf_x=\l f_x,\qqq \
(\l,x)\in \mathbb{C}\ts\Z_+,
\]
with initial conditions
 \[\varphi_0 = \vartheta_1 = 0, \;\;\; \varphi_1 = \vartheta_0 = 1.
\]
Note that for $z\in \mathbb{D}$ the functions $\varphi_x(\l(z))$
and $\vartheta_x(\l(z))$ are polynomials in $z+z^{-1}$. Therefore,
they can be considered outside of $\mathbb{D}$ too. Now we can write
down the Wronskian for the Jost and fundamental solutions. For $x =
0$ we have
$$
\{\psi, \varphi\}|_{x=0} = \psi_0\varphi_1 - \psi_1\varphi_0= \psi_0.
$$
Since the Wronskian of two fundamental solutions does not depend on $x,$ we have at $x = p:$
$$
\psi_0 = \{\psi, \varphi\}|_{x=p} = a_p(\psi_p\varphi_{p+1} -
 \psi_{p+1}\varphi_p) = z^p\varphi_{p+1} - a_pz^{p+1}\varphi_p.
 $$
Similar arguments for $\psi$ and $\vt$ yeild
\[\lb{wro1}
\psi_1(z) = -z^p\vartheta_{p+1}(\l(z)) + a_pz^{p+1}\vartheta_p(\l(z)).
\]
Recall that $\psi^{\diamond}(z) = \psi(\frac{1}{z}).$  Computing the
Wronskian (\ref{wro}) for $\p$ and $\psi^{\diamond}$ and using the
fact that $\varphi(\l(z)) = \varphi(\l(\frac{1}{z})),$ we obtain
\[\lb{fp1}
\textstyle \varphi_{p+1}(\l(z))(z^2 - 1) = z^{p+2}\psi_0
\Bigl(\frac{1}{z}\Bigr) - \frac{\psi_0(z)}{z^p},
\]
\[\lb{fp2}
\textstyle
 a_p\varphi_{p}(\l(z))(z^2 - 1) =
 z^{p+1}\psi_0\Bigl(\frac{1}{z}\Bigr) - \frac{\psi_0(z)}{z^{p-1}}.\]
Recall that

$(\m_j)_0^{p-1}$ are the roots of $\vp_{p+1}(\l)$ and  the
eigenvalues of $J_p,$

$(\t_j)_1^{p-1}$ are the roots of $\vp_p(\l)$ and the  eigenvalues
of $J_{p,1},$

$(\vr_j)_1^{p-1}$ are the roots of $\vt_{p+1}(\l)$ and  the
eigenvalues of $J^1_p,$

$(\n_j)_1^{p-2}$ are the roots of $\vt_p(\l)$ and the eigenvalues of $J^1_{p,1}.$

We define two functions $u^{\pm}$ by
\[\lb{uv} u^{\pm} = \varphi_{p+1} \mp a_p \varphi_p.
\]
Denote the characteristic polynomials of the matrices  $J_p,
J_{p,1}$ and $J^{\pm}$ by $\mathcal{D}_p, \mathcal{D}_{p,1}$ and
$\mathcal{D}^{\pm}$ respectively. In order to prove Theorems
\ref{pr1} and \ref{T4}, we recall a well-known result about
fundamental solutions and study the connection between functions
$u^{\pm}(\l)$ and matrices $J^{\pm}.$

\begin{lemma}\lb{ldop}
Let $J$ be a Jacobi operator acting on $\ell^2(\mathbb{N})$  defined
by (\ref{j}) with a perturbation $q(a,b) \in \gX_{k},\; k\in \{2p-1,
2p\}.$ Then

\no i) the fundamental solutions $\vp(\l)$ satisfy
\[\lb{fi1}
\textstyle
 \varphi_x(\l) = \frac{\l^{x-1}}{a_1 \cdots
a_{x-1}}\Bigl(1- \frac{b_1+\dots+b_{x-1}}{\l}+ O(\l^{-2})\Bigr)\;\;
\text{ as } \l\rightarrow+\infty, \;\; x\in \mathbb{N}_p,
\]
ii) the functions $u^{\pm}(\l)$ defined by (\ref{uv})  satisfy
$\mathcal{D}^{\pm} = (-1)^pA_1a_p u^{\pm}.$
\end{lemma}
\textbf{Proof.} i) Using  the equation
$a_{x}f_{x+1}=(\l-b_x)f_x-a_{x-1}f_{x-1}$ and $\vp_0 = 0,\; \vp_1 =
1$, we obtain: \no if $x=1$, then we have
\[
\lb{asp0} \textstyle \vp_2={(\l-b_{1})\vp_{1}-\vp_{0}\/a_1} ={\l -
b_1\/a_1},
\]
if $x=2$, then we have
\[
\lb{asp0} \textstyle \vp_{3}={(\l-b_{2})\vp_{2}-a_1\vp_{1}\/a_{2}}
={(\l-b_{2})(\l - b_1)-a_1^2\/a_1a_2} ={\l^2-(b_1+b_2)\l+b_1b_2 -
a_1^2\/a_1a_{2}}.
\]
Repeating this procedure we obtain \er{fi1}.

ii) Recall that the matrices $J^{\pm}$ are given by
\[\lb{jpm}J^{\pm} = \left(\begin{array}{ccccccccc}
b_1 & a_1 & 0 & 0 &  ...  &  0\cr
 a_1 & b_2 & a_2 & 0 &  ... &  0 \cr
 0 & a_2 & b_3 & a_3 &  ... & 0 \cr
 ... & ... & ... & ... & ... & ...  \cr
 0 & ...  & 0 & a_{p-2} & b_{p-1} & a_{p-1}\cr
 0 & ...  & 0 & 0 & a_{p-1} & b_p\pm a_p^2 \cr
\end{array}\right) .
\]
It follows from the asymptotics (\ref{fi1}) that
\[
\textstyle
 \mathcal{D}_p = (-1)^{p}A_1\vp_{p+1}, \qqq
\mathcal{D}_{p,1} = \frac{(-1)^{p-1}A_1}{a_p}\vp_{p}.
\]
Using the linearity of the determinant over rows, we divide the last
row of $J^{\pm}$ into two parts and obtain
$$
\mathcal{D}^{\pm} = \mathcal{D}_p \pm a_p^2\mathcal{D}_{p,1}
= (-1)^pA_1a_p(\vp_{p+1} \mp a_p\vp_{p}) = (-1)^pA_1a_p u^{\pm}.
$$
\BBox

Recall a corollary from \cite{HJ12}, see page 241.

\begin{lemma}\lb{book}
Let $n>2,$ let $A$ be Hermitian $n\times n$ matrix and  let $v\in
\mathbb{C}^n$ be nonzero. If $A$ and $A+vv^*$ do not have similar
eigenvectors, then
$$
\l_j(A)< \l_j(A+vv^*)< \l_{j+1}(A), \;\; j\in\mathbb{N}_{n-1},
$$
$$
\l_n(A)< \l_n(A + vv^*).
$$
\end{lemma}
\begin{lemma}Let $J$ be a Jacobi operator acting
on $\ell^2(\mathbb{N})$ defined by (\ref{j}) with a perturbation
$q(a,b) \in \gX_{2p}, \;p\in \mathbb{N}.$ Let the operators $J^{\pm}$
be defined by (\ref{P1}) and let $(\a_j^{\pm})^p_1$ be their
eigenvalues. Then the relations (\ref{ord}-\ref{ord2}) hold true and
\[ \begin{array}{cc}
    \alpha^+_p \rightarrow +\infty, \;\; \alpha^-_1
     \rightarrow -\infty & \;\;\;\;\;\text{as }\;\; a_p
      \rightarrow +\infty, \\
    \alpha^{\pm}_n \rightarrow \mu_{n-1} \;\;\forall
    n\in \mathbb{N}_{p} & \text{as }\;\; a_p \rightarrow 0.
  \end{array}
\]
\end{lemma}
\no \textbf{Proof.} The relations (\ref{ord}-\ref{ord2})  follow
from Lemma \ref{book} if we take $v = (0,\dots, 0, a_p).$ Note that
the fundamental solution $\varphi_p(\l)$ does not depend on $a_p,$
while $\varphi_{p+1} = \frac{\varphi_{p+1}^o}{a_p},$ where
$\varphi_{p+1}^o$ is the corresponding fundamental solution for $a_p
= 1.$ We will consider the eigenvalues $\alpha^+_n,$ the proof for
$\alpha^+_n$ is similar.

i) Let $a_p \rightarrow +\infty.$ It follows from (\ref{fi1})  that
the second term of the polynomial $\varphi_{p+1}(\l)$ is equal to
$-\frac{B_0}{A_1}\l^{p-1},$ and the first term of
$-a_p\varphi_{p}(\l)$ is equal to $-\frac{a^2_p}{A_1}\l^{p-1}.$ The
sum of these terms is the second term of $u^+(\l),$ and after
dividing it by the leading term of $u^+(\l)$, which is equal to the
leading term of $\vp_{p+1}(\l)$, its absolute value tends to
infinity as $a_p$ tends to infinity. Then the absolute value of the
sum $\sum_1^p \alpha^+_j$ of its roots also tends to infinity. Since
all roots of $u(\l),$ except for $\alpha^+_p,$ are bounded, we have
$\alpha^+_p \rightarrow +\infty.$

ii) Let $a_p \rightarrow 0.$ We can see from (\ref{jpm})  that the
matrices $J^{\pm}$ are obtained from the matrix $J_p$ through the
one-dimensional perturbation. When $a_p$ tends to 0, this
perturbation also tends to 0. Thus, the eigenvalues $\alpha^{\pm}_n$
of these matrices tend to those of $J.$ \BBox

In order to prove the main theorems, we need the following  property
of the phase shift $\x(z).$
\begin{lemma}\lb{lz}
Let the operator $J$ be defined by (\ref{j})  with a perturbation $q
\in \gX^+_{k},$ where $k\geq 1$,  and let the function $\x(z)$ be
its phase shift. Then the relations (\ref{lz1}-\ref{lz2}) hold true.
\end{lemma}
\no {\bf Proof.} We use arguments from \cite{APR19}, where similar
relations are proved in  for discrete Schrodinger operator. We know
that $\psi_0$ has $\gN$ simple  zeros in $\mathbb{D},$ all of them
on the real line, and possibly a simple zero in $\pm 1.$ Applying
the argument principle to $\psi_0(z)$ along the unit circle, we see
that the total variation of $\arg\psi_0$ along the unit circle in
the counterclockwise direction is given by
\[\lb{lzd1}
\text{var} \arg \psi_0 = 2\pi\Bigl(\gN +\frac{m_++m_-}{2}\Bigr).
\]
Since for $|z| = 1$ we have $\x(\frac{1}{z}) = \x(\bar{z}) =
-\x(z),$ half of this change is obtained  as $z$ moves from 1 to $-1$.
Therefore, the total variation of $\arg \psi_0$ along the unit
circle from 1 to $-1$ in the counterclockwise direction is given by
\[\lb{lzd2}
\text{var}_+ \arg \psi_0 = \pi\Bigl(\gN +\frac{m_++m_-}{2}\Bigr).
\]
One can see that $\psi_0(1)$ is real, and the sign of it  depends on
the amount of zeros in $(0,1].$ Therefore, we have $\x(1+0i) = -\pi
n_++\frac{\pi m_+}{2}.$ Adding the change $(\ref{lzd1}-\ref{lzd2}),$
we obtain the relations (\ref{lz1}-\ref{lz2}). \BBox

\section {Proof of main theorems}
\setcounter{equation}{0} \textbf{Proof of Theorem \ref{T1}.} It
follows from  (\ref{t12}) that if $\alpha^+_p\leq 2$, then
$\l_{n_+}<\alpha^+_p,$ and there are no eigenvalues in $[2,
+\infty).$ If $\alpha^+_p = 2,$ then there is a virtual state $z =
1,$ but still no positive eigenvalues. We fix a perturbation
$q_0(a_0, b_0)$ from $\mathfrak{X}_{2p}$ or $\mathfrak{X}_{2p-1}$,
such that $\alpha^+_p(q_0) = 2.$ Define the perturbation
$q_{\varepsilon} = q_{\varepsilon}(a_{\varepsilon},
b_{\varepsilon}), \; \varepsilon\in \mathbb{R}$ by
\[(q_{\varepsilon})_x = 0,\qq x>2p, \qqq (a_{\varepsilon},
b_{\varepsilon})_x = (a_x, b_x + \varepsilon),\qq x\leq p.
\]
We sometimes write $\psi(z, \varepsilon), \alpha^+_n(\varepsilon),$
etc. instead of $\psi(z), \alpha^+_n,$ etc. when several
perturbations $q_{\varepsilon}$ are being dealt with. The operator
corresponding to the perturbation $q_{\varepsilon}$ we denote by
$J(\varepsilon).$ Note that the eigenvalues
$\alpha^+_n(\varepsilon), \alpha^-_n(\varepsilon)$ satisfy
\[\lb{eps}\alpha^+_n(\varepsilon) = \alpha^+_n(0)+\varepsilon,
\qqq \alpha^-_n(\varepsilon) = \alpha^-_n(0)+\varepsilon.
\]
Due to (\ref{wro}) the Jost function for $J(\varepsilon)$ for
$\lambda = z+ \frac{1}{z},\;\lambda\in(0, +\infty)$ has the form
\[\lb{wre}
\psi(z, \varepsilon) = z^p\varphi_{p+1}(\l,\varepsilon)
 - a_pz^{p+1}\varphi_p(\l,\varepsilon) = z^p\varphi_{p+1}
 (\l-\varepsilon,0)- a_pz^{p+1}\varphi_p(\l-\varepsilon,0).
 \]
The function $\psi(z, \varepsilon)$ is entire in $z,\varepsilon$
and satisfies
\[\psi(1,0) = 0,\qq \psi_t(1,0)\neq0.\]
Then due to Implicit Function Theorem there exists a function $z(\varepsilon)$, analytic in small disk $\{|\varepsilon|<\epsilon\}$ such that $\psi(z(\varepsilon),\varepsilon) = 0$ in the disc $\{|\varepsilon|<\epsilon\}$. Here $\l_{n_+}(\varepsilon) = z(\varepsilon)+\frac{1}{z(\varepsilon)}$ for $\varepsilon>0.$

We have $z(\varepsilon) = 1+z_1\varepsilon+O(\varepsilon^2).$  Then
from (\ref{wre}) we obtain for $z = z(\varepsilon)$ as
$\varepsilon\rightarrow0:$
\[\varphi_{p+1}(\l-\varepsilon) = a_p z(\varepsilon)
\varphi_{p}(\l-\varepsilon) = a_p(\varphi_{p}(2) +
\varepsilon(z_1\varphi_p(2)-\dot{\varphi}_{p}(2)))+
O(\varepsilon^2),
\]
\[\varphi_{p+1}(\l-\varepsilon) = \varphi_{p+1}(2) -
 \varepsilon\dot{\varphi}_{p+1}(2)+O(\varepsilon^2),
\]
which yields $z_1 = -\frac{\dot{\varphi}_{p+1}(2)
- a_p\dot{\varphi}_p(2)}{a_p\varphi_p(2)}<0,$
where $\dot{u} = \frac{d}{d\lambda}u.$
If $\varepsilon>0$ is small enough, then (\ref{eps})
yields $\alpha^+_p(\varepsilon) = \varepsilon$ and
 $\alpha^-_p(\varepsilon) = \alpha^-_p(q_0)+\varepsilon<2.$
Thus we obtain $\alpha^-_p(\varepsilon)<2<\lambda_{n+}<\alpha
_p(\varepsilon).$ If $\varepsilon$ is increasing, then all
eigenvalues $\alpha^-_p(\varepsilon)<2<\lambda_{n+}<\alpha
_p(\varepsilon)$ move monotonically to the right.

We define a matrix $J^1\in\gX_{k-1}$ which is be obtained from $J$
by deleting the first row and the first column. One can check that
the Jost function for $J^1$ is equal to $\frac{\p_1(z)}{z}.$ It is
expressed in terms of the fundamental solutions by
\[\lb{wrrr}
\frac{\psi_1(z)}{z} = a_pz^p\vartheta_p(\l(z)) - z^{p-1}
\vartheta_{p+1}(\l(z)).
\]
Denote the eigenvalues of the matrices $(J^1)^{\pm}$ as
$(\gamma^{\pm})_1^{\lfloor k \rfloor}$ and denote the eigenvalues of
$J^1$ as $\omega_j.$

At $\varepsilon_1 = 2-\gamma^+_{p-1}(0)$  we have
$\gamma^+_{p-1}(\varepsilon) = 2.$

Note that if $\gamma^+_{p-1}\leq 2,$ then from (\ref{wrrr})  we
deduce that the operator $J^1$ does not have eigenvalues in $(2,
+\infty).$ Let $\varepsilon = \varepsilon_1+y,$ where $\varepsilon_1
= 2-\gamma^+_{p-1}(0).$ Then due to \ref{wrrr} we obtain for
$\lambda = z+ \frac{1}{z},\lambda\in(0, +\infty):$
\[\lb{wre1}
F(z,y) = \frac{\psi_1(z,\varepsilon)}{z}
= a_pz^p\vartheta_p(\l-y,\varepsilon_1) -
z^{p-1}\vartheta_{p+1}(\l-y,\varepsilon_1).
\]
The function $F(z,y)$ is entire in $z,y$ and satisfies
\[
F(1,0) = 0, \;\;F_z(0,0)\neq 0.
\]
Due to Implicit Function Theorem there is  a function $z_o(y),$
analytic in small disc $\{|y|<\epsilon_o\}$ such that $F(z_o(y),y) =
0$ in the disc $\{|y|<\epsilon_o\}.$ Here we have
$\omega_{m_+}(\varepsilon) = z_o(y)+\frac{1}{z_o(y)}, y =
\varepsilon - \varepsilon_1.$ We have $z_o(y) = 1+z_1y+O(y^2)$ as
$y\rightarrow0.$ Then from (\ref{wre1}) we obtain for $z = z_o(y)$
as $y\rightarrow0:$
\[
\vartheta_{p+1}(\l-y) = a_p z_o(y)\vartheta_{p}(\l-y,\varepsilon_1)
= a_p(\vartheta_{p}(2,\varepsilon_1) + y(z_1\vartheta_p(2,\varepsilon_1)
-\dot{\vartheta}_{p}(2,\varepsilon_1)))+O(y^2),
\]
\[\vartheta_{p+1}(\l-y) = \vartheta_{p+1}(2,\varepsilon_1) - y\dot{\vartheta}_{p+1}(2,\varepsilon_1)+O(y^2),
\]
which yields
$z_1 = \frac{-\dot{\vartheta}_{p+1}(2,\varepsilon_1)
 + a_p\dot{\vartheta}_p(2,\varepsilon_1)}
 {a_p\vartheta_p(2,\varepsilon_1)}<0.$

If $y>0$ is small enough, then $\gamma^+_{p-1}(\varepsilon)
 = 2+y$ and $\gamma^-_{p-1}(\varepsilon) = \gamma^-_{p-1}(0)
+\varepsilon<2,$ since there is the basic relation (\ref{ord}). Thus
we obtain
$$
\gamma^-_{p-1}(\varepsilon)<2<\omega_{m_+}(\varepsilon)<\gamma^+_{p-1}
(\varepsilon), \qq \omega_{m_+}(\varepsilon)<\lambda_{n_+}(\varepsilon).
$$
It is important that $\omega_{m_+}(\varepsilon)<\lambda_{n_+}
(\varepsilon)$ for any $\varepsilon$, since
$\sigma_d(J)\bigcap\sigma_d(J^1) = \emptyset.$ If $\varepsilon$ is
increasing,  the eigenvalues
$\gamma^-_{p-1}(\varepsilon)<\omega_{m_+}(\varepsilon)
<\gamma^+_{p-1}(\varepsilon)$ and
$\omega_{m_+}(\varepsilon)<\lambda_{n_+}(\varepsilon)$ move
monotonically to the right. New eigenvalues do not appear until
$\varepsilon = \varepsilon_2, \;\alpha^+_{p-1}(\varepsilon_2) = 0.$ We
can repeat the above arguments for any $q_0$ and $\varepsilon$ to
get the other inequalities. The proof for the interval $(-\infty,
-2)$ is similar. Thus, we obtain (\ref{t11}). The estimate
(\ref{t12}) follows from Lemma \ref{f02}. \BBox

In order to prove Proposition \ref{TB}, we need the following
Hochstadt result from \cite{H78}.

\begin{theorem}\lb{Ta}
Let $J_p$ be a finite Jacobi matrix on $\C^{p}$ defined  by
(\ref{J1}). Let $(\m_j)_0^{p-1}$ and $(\t_j)_1^{p-1}$ be its mixed
and Dirichlet eigenvalues. Then the mapping $\mu\star\tau:
\mathbb{R}_+^{p-1} \times \mathbb{R}^{p} \rightarrow \cE_{2p-1}$
given by $q = (a,b)\mapsto \mu\star\tau$ is a real analytic
isomorphism between $\mathbb{R}_+^{p-1} \times \mathbb{R}^{p}$ and
$\cE_{2p-1}$.

Moreover, there exists an algorithm  to recover the Jacobi matrix
from its eigenvalues.\end{theorem}

\textbf{Remark.} Analyticity of the mapping was proved by Simon  and
Gesztesy \cite{GS97}. The recovery algorithm is given below, in the
proof of Proposition \ref{TB}.

\textbf{Proof of Proposition \ref{TB}}. It follows from Lemma
\ref{ldop} that the numbers $(\alpha^{\pm}_j)_1^p$ are the roots of
the polynomials $u^{\pm}(\l)$ and that we have
\[\lb{TB3}
\mathcal{D}_p = \frac{\mathcal{D}^+ +
 \mathcal{D}^-}{2}, \;\;\; \mathcal{D}_{p,1} =
 \frac{\mathcal{D}^+ - \mathcal{D}^-}{2a_p^2}.
 \]
Also note that $J_{p,1} = J^{\pm}_1.$ Now we determine the coefficients of $J.$ Using the trace formula, we obtain
\[\lb{TB4}
\Tr J^+ - \Tr J^- = \Sigma_1^p \alpha^+_j - \Sigma_1^p
\alpha^-_j = \big(\Sigma_1^p b_j + a_p^2\big)-
\big(\Sigma_1^p b_j - a_p^2\big) = 2a^2_p.
\]
Since the coefficient $a_p$ is positive, it is determined. It follows from the trace formulas for the matrices $J_p$ and $J_{p,1}$ that the second term of $\chi(\l)$ equals $-\Sigma_1^p b_j,$ while the one of $\chi_1(\l)$ equals $-\Sigma_1^{p-1} b_j.$ Taking the difference, we obtain $b_p.$ All the other coefficients are the same for the matrices $J_p$ and $J^+.$ Thus, we can find the coefficients of the matrix $J^+$ and therefore find the matrix $J_p$ too. We use a standart algorithm from \cite{H78}. Let the last $k$ rows of the matrix $J^+$ be already found. Denote the standart basis in $\mathbb{C}$ as $(\delta_j)_1^p.$ From \cite{H78} we obtain that
\[\lb{Tb1}
((J^+)^k \delta_p, \delta_p) =-
\sum_1^p \frac {(\alpha_j^+)^k \mathcal{D}_{p,1}(\alpha^+_j)}
{(\mathcal{D}^+)'(\alpha^+_j)},\;\; k = 1, \dots, p.
\]
On the other hand, direct calculations show that
\[\lb{Tb2}
((J^+)^{2k+1} \delta_p, \delta_p) = d_1 b_{p-k}+d_2 a_{p-k}+\dots d_{2k-2} a_{p-1} +d_{2k-1} b_p,
\]
where the coefficients $d_j$ depend only  on $a_{p-k}, b_{p-k+1},
\dots, a_{p-1}, b_p$ and therefore are known, and $d_1\neq0.$ Thus,
the right parts of (\ref{Tb1}) and (\ref{Tb2}) are equal, and
everything except $b_{p-k}$ is known, which allows us to determine
$b_{p-k}$ from the obtained equation. After that we can consider
$((J^+)^{2k+2} \delta_p, \delta_p)$ and use a similar principle to
find $a_{p-k-1},$ so now the last $k+1$ rows are known. Proceeding
like this, we can find all the elements of the matrix $J^+$ and,
since we already found $a_p$ and $b_p,$ of the matrix $J.$ The
continuity of the corresponding mapping follows from the results of
\cite{H78}. \BBox

In order to prove Theorem \ref{pr1}, we need  to determine the
derivative of the phase shift $\x(z).$
\begin{lemma}\lb{lpr}Let the operator $J$ be defined by
(\ref{j}) with a perturbation $q \in \gX^+_{k},$ where $k\in \{2p-1,
2p\}$. Then for $\l\in(-2,2)$ its phase shift function $\x(z)$
satisfies
 \[\lb{ss}
 \textstyle
\sqrt{1 - {\l^2\/4}} (\xi(z(\l)))' =
\sum_{r_j\in \mathbb{R}}\frac{r_j \Re z - 1}{2|z - r_j|^2}+
\sum_{r_j\in\mathbb{C}^+}\frac{\Re r_j \Re z - 1}{|z - r_j|^2},
\]
where $z(\l) = {\l\/2}-\sqrt{{\l^2\/4}-1}$ and $(r_j)_1^{k}$  are
the roots of the Jost function $\psi_0(z).$
\end{lemma}
\no \textbf{Proof.} Consider the S-matrix
\[\lb{s}
S(z) = \frac{\psi_0(z^{-1})}{\psi_0(z)} = e^{-2i\xi(z)}.
\]
where $|z| = 1,\; \Im z <0.$
Taking the derivative of the both parts of (\ref{s}), we obtain:
\[
-\frac{\psi_0'(\frac{1}{z})z'(\l)}{z^2\psi_0(z)}
- \frac{\psi_0(\frac{1}{z})\psi_0'(z)z'(\l)}{\psi_0^2(z)}
= -2ie^{-2i\xi(z)}\xi'(z)z'(\l),
\]
which yields
\[
\frac{\psi_0'(\frac{1}{z})}{z^2\psi_0(\frac{1}{z})}  +
\frac{\psi_0'(z)}{\psi_0(z)} = 2i\xi'(z).
\]
From Lemma \ref{L1} we have $\psi_0(z) = C\Pi_{j = 1}^k (z - r_j).$ Then
\[
\frac{\psi_0'(z)}{\psi_0(z)} = \sum_{j = 1}^k\frac{1}{z - r_j}
 =\sum_{j = 1}^k\frac{\bar{z} - \bar{r_j}}{|z - r_j|^2},
 \]
and
\[
\frac{\psi_0'(\frac{1}{z})}{z^2\psi_0(\frac{1}{z})}
+ \frac{\psi_0'(z)}{\psi_0(z)} =
\sum_{j = 1}^k\frac{\bar{z} - \bar{r_j}+
 \bar{z}^2(z - r_j)}{|z - r_j|^2} =
 \sum_{j = 1}^k\frac{2\bar{z} - \bar{r_j}
 (1 + \bar{z}^2)}{|z - r_j|^2}.
\]
Recall that $z'(\l)=-{z\/2\sqrt{{\l^2\/4}-1}}$  (here we use the
positive branch of the square root). Then for the full derivative we
have
\[\lb{s1}
\sqrt{1 - {\l^2\/4}} \xi'(z)z'(\l) =
-\sum_{j = 1}^k\frac{z(2\bar{z} -
\bar{r_j}(1 + \bar{z}^2))}{4|z - r_j|^2} =
\sum_{j = 1}^k\frac{\bar{r_j} \Re z - 1}{2|z - r_j|^2}.
\]
Since $\psi(z)$ is a polynomial with real coefficients, for every complex root $r_j$ there is also a root $\bar{r_j}.$ Taking the sum of the corresponding parts of \ref{s1}, we obtain \ref{ss}.  \BBox

\textbf{Proof of Theorem \ref{pr1}.}
Substitution of (\ref{fp1}-\ref{fp2}) into (\ref{uv}) gives
\[\lb{tpr11}
(z^2-1)u^+(z(\l)) =
(z-1)(z^{p+1}\p_0\Bigl(\frac{1}{z}\Bigr)+\frac{1}{z^p}\p_0(z)),
\;\;\l\in(-2,2).
\]
Recall that the roots of $u^+(\l(z))$ are $\a^+_j, \; j\in\mathbb{N}_{p}.$ For $|z|=1, \;\Im z< 0$ we rewrite $(z+1)u^+(\l(z))=0$ by
\[\lb{pr12}
\psi_0(z)(z^{2p+1}S(z) + \bar{z})=0.
\]
Since we consider $|z|=1,$ we can express (\ref{pr12})  in terms of
arguments. This, we obtain that $z(\a^{+}_j)$ are the roots of
\[\frac{2\xi(z)}{2p+1} +\arg z = -\frac{(2\pi-1) n}{2p+1}, \;\; n
\in \N_p.\] The proof for $\alpha^-_j$ is similar. Let $F(\l) =
\frac{2\xi(z(\l+0i))}{2p+1} +\arg z(\l+0i).$ Using the result of
Lemma \ref{lpr}, we obtain that
$$
\sqrt{1 - {\l^2\/4}}z'(\l+0i)F'(\l) =\frac{1}{2}+
 \frac{1}{2p+1}\sum_{j = 1}^k\frac{\bar{r_j}
\Re z - 1}{2|z - r_j|^2} \geq \frac{1}{4p+2}\sum_{j =
1}^k\Bigl(\frac{\bar{r_j} \Re z - 1}{|z - r_j|^2} + 1\Bigr).
$$
For real roots $r_j$ and $|z|=1$ we have
\[\frac{\bar{r_j}
\Re z - 1}{|z - r_j|^2} + 1 = \frac{|z|^2 + |r_j|^2 -
r_j\Re z - 1}{|z - r_j|^2} = \frac{|r_j|^2 -r_j\Re z}{|z - r_j|^2} >0,
\]
since $|r_j|\geq1$. The proof for non-real roots is  similar.
Therefore, $F'(\l)>0$ for $\l \in (-2, 2),$ so the function $F(\l)$
is strongly increasing. We also have $F(-2) = -\pi$ and $F(2) = 0,$
so all of the roots of $F(\l)+\frac{(2\pi-1) n}{2p+1}$ for
$n\in\mathbb{N}_p$ lie in $[-2, 2].$ That means that all $p$ roots
of $u^+(z)$ are the roots of $F(\l)+\frac{(2\pi-1) n}{2p+1}$ for
such $n.$ The proof for $u^-(z)$ is similar. It suffices to use
Proposition \ref{TB} and obtain a bijection. \;\BBox

In order to prove Theorem \ref{T4}, we need a simple corollary of
Theorem \ref{T1}.

\begin{corollary}\lb{lemm}
Let $J$ be a Jacobi matrix given by (\ref{j}) with a perturbation $q
\in\gX_{k}$ and let $(s_j), \; j\in n_{\bu}$ be the roots of its
Jost function in  $\mathbb{D}.$ Let $(\mu_n)_0^{p-1}$ be eigenvalues
of the matrix $J_p$.

i) If $\m_{p-1}>2$, then we have $s_1 < z(\m_{p-1}).$

ii) If $\m_{0}<-2$, then we have  $s_{-1} >
z(\m_{0}).$

\end{corollary}

\no {\bf Proof.} If $\m_{p-1}>2,$ then \er{t12} gives that $\l_{1}>
\mu_{p-1}.$ Since the function $z(\l)$ decreases when $\l>2,$ we
have $s_1 = z(\l_{1}) < z(\m_{p-1}),$ which proves i). Proof of ii)
is similar. \BBox

\textbf{Proof of Theorem \ref{T4}}. First of all, we need to show
that for every  operator $J$ with perturbation $q\in \gX_{k}$ we
have $(r_j)_1^k\in \mathcal{R}_k.$ It follows directly from Lemma
2.2.

After that we prove the injection. Let the set $(r_n)_1^{k}$  be
given. Then the polynomial $\psi_0(z)$, which roots are exactly
$(r_n)_1^{k},$ counting multiplicities, is determined up to
multiplication by a constant. We can use Lemma \ref{f01} to find
$a_p$. Since we have (\ref{fp1}-\ref{fp2}), the functions
$\varphi_{p+1}(z)$ and $\varphi_{p}(z)$ are also determined up to
multiplication by a constant, which means that their roots are
uniquely determined. It suffices to use the fact that these roots
define the matrix $\cJ_p$ uniquely, which follows from Theorem
\ref{Ta}.

Now we prove the surjection. We start with the case  when
$|\psi_0(z)|>0$ in $\dD.$ It follows from Theorem \ref{pr1} that
there is a matrix $J$ which Jost function is exactly $\psi_0(z).$ We
can recover it using the algorithm from Theorem \ref{pr1}.

Now let us use induction by the pair of numbers $k$ and the  amount
$m$ of roots of $\psi(z)$ in $\dD.$ The case $k = 0$ is already
proved, the case $k = 1$ is obvious. Fix $k$ and $m$ and assume that
the surjection is proved for any $(\tilde{m},\tilde{k})< (m,k).$ Let
$k = 2p$ be even, the proof for the case $k = 2p-1$ is similar.
Consider a function $\psi(z)$ with $m$ zeros in $\dD_1$ and let
$s^o$ be one of its zeros with the smallest absolute value. Consider
the function $\psi_0^o(z) = \frac{\psi_0(z)}{z - s^o.}$ It is easy
to see that $\psi_0^o(z)$ also suffices R1-R3 and that it has one
less root in $\dD_1.$ It means that there is a finite Jacobi matrix
$J_p^o$ such that $\psi_0^o(z)$ is the Jost function of the
corresponding matrix $J^o$. Denote the fundamental solutions for
this matrix as $\varphi^o_{p+1}(\l)$ and $\varphi^o_{p}(\l).$
Substituting $\psi_0^o$ into (\ref{fp1}), we obtain
\[
\varphi_{p+1}(\l) = -s^o \varphi^o_{p+1}(\l) + a_p\varphi^o_{p}(\l).
\]
Denote the roots of $\varphi^o_{p+1}(\l)$ as $(\mu^o_n)_0^{p-1}.$ They are all real and simple. Then we have
\[
\varphi_{p+1}(\mu^o_n) = a_p\varphi^o_{p}(\mu^o_n).
\]
Since the roots of $\varphi^o_{p+1}(\l)$ and $\varphi^o_{p}(\l)$
 alternate, we have
 \[
\varphi^o_{p}(\mu^o_{n-1})\varphi^o_{p}(\mu^o_{n}) < 0,\;\; n =
0,1,\dots,p-1,
\]
and, consequently,
\[
\varphi_{p+1}(\mu^o_{n-1})\varphi_{p+1}(\m^o_{n}) < 0,\;\;
n = 0,1,\dots,p-1.
\]
This means that there is a root of $\varphi_{p+1}$ between  any two
roots of $\varphi^o_{p+1}.$ Since the sign of $\varphi_{p+1}(\l)$ at
$\l \rightarrow -\infty$ is similar to the one of
$\varphi^o_{p+1}(\l)$ and different from the sign of
$\varphi^o_{p}(\m^o_0)$, the last zero of $\varphi_{p+1}$ has to be
smaller that $\m^o_0.$ Thus, all zeros of $\varphi_{p+1}$ are real
and simple. Similar arguments can be applied to obtain that
\[\lb{m0}
a_p\varphi_{p}(\l) = a_p(\l - s^o)\varphi^o_{p}(\l) - \varphi^o_{p+1}(\l).
\]
and that all roots of $\varphi_{p}(\l)$ are real and simple. It
suffices to show that the roots of $\varphi_{p+1}(\l)$ and
$\varphi_{p}(\l)$ alternate. Let $\mu$ be a root of
$\varphi_{p+1}(\l).$ Then we have
\[\lb{m1}
a_p\varphi^o_{p}(\m) = s^o \varphi^o_{p+1}(\m).
\]
Substituting (\ref{m1}) into (\ref{m0}), we obtain that
\[\lb{vpp}
\varphi_{p}(\mu) = \varphi^o_{p+1}(\mu)(s^o\m - (s^o)^2 - 1).
\]
It follows from Lemma \ref{lemm}  that $\m^o_{p-1}< s^o +
\frac{1}{s^o },$ and we have already proved that $\m^o_{p-1}>
\m_{p-1}.$ Then $\m_{p-1}< s^o + \frac{1}{s^o },$ so the sign of the
second factor in (\ref{vpp}) is the same for all eigenvalues $\m.$
Then we have \[\lb{alt}\sign\varphi_{p}(\mu) = \sign \varphi^o_{p+1}(\m).\]
We already know that the roots of $\varphi_{p+1}$ and
$\varphi^o_{p+1}$ alternate, so \ref{alt} implies that the roots of
$\varphi_{p}$ and $\varphi_{p+1}$ also alternate. Thus, it follows
from Theorem \ref{Ta} that there is a finite Jacobi matrix $J_p$
such that $\varphi_{p}$ and $\varphi_{p+1}$ are its fundamental
solutions.  We can use Lemma \ref{f01} to find $a_p.$ Adding the
element $a_p$ to the matrix $J_p,$ we obtain the matrix $J,$ such
that $(r_n)_1^{k}$ is the sequence of roots of its Jost function.
Therefore, the surjection is proved.

Therefore, the recovery algorithm looks like follows:

1. Determine  $a_p$ using Lemma \ref{f01}.

2. Using the sequence $(r_n)_1^{k},$ find the function  $\psi_0(z)$
up to multiplication by a constant.

3. Determine  the functions $\vp_{p}(z)$ and $\vp_{p+1}(z)$ up
to multiplication by a constant using the formulas
(\ref{fp1}-\ref{fp2}).

4. Calculate the eigenvalues $(\m_j)_0^{p-1}$ and $(\t)_1^{p-1}$.

5. Recover the matrix $J_p$ using Theorem \ref{Ta}.
\BBox

\textbf{Proof of Corollary \ref{corxi}.} We have $\frac{\p_0(q_1,
w_j^{-1})}{\p_0(q_1, w_j)}  = \frac{\p_0(q_2, w_j^{-1})}{\p_0(q_2,
w_j)}.$ Therefore, we can obtain that
\[
w_j^{k}\p_0(q_1, w_j^{-1})\p_0(q_2, w_j) =
 w_j^{k}\p_0(q_2, w_j^{-1})\p_0(q_1, w_j).
\]
Note that these equalities also hold true if we swap $w_j$  and
$w^{-1}_j.$ Consider the functions $G_1(z) =z^{k}\p_0(q_1,
z^{-1})\p_0(q_2, z)$ and $G_2(z) = z^{k}\p_0(q_2, z^{-1})\p_0(q_1,
z).$ These are polynomials of degree $2k$ which coincide at $2k+1$
points (including 1). Therefore, these polynomials are equal.
Consider the identity
\[\lb{korni}z^{k}\p_0(q_1, z^{-1})\p_0(q_2, z)
 = z^{k}\p_0(q_2, z^{-1})\p_0(q_1, z).
 \]
Let $r$ be some zero of $\p_0(q_2, z).$ It  is also a zero of the
right part of (\ref{korni}). It follows from Lemma \ref{f02} that
$r\neq 0$ and $\p_0(q_2, r^{-1})\neq 0.$ Therefore, we have
$\p_0(q_1, r) = 0.$ We can swap $\p_0(q_1, z)$ and $\p_0(q_2, z)$
and obtain that these polynomials have similar roots. It suffices to
use Theorem \ref{T4} and obtain that $J_1 = J_2.$
\BBox

\section {Location of resonances}
\setcounter{equation}{0} Here we state and proof the results about
location of resonances.  We give estimates on absolute value of
resonances of $J$ depending on its coefficients and the eigenvalues
of $J$ and $J^{\pm}.$ We also consider the case when one of the
resonances of $J$ is moved to infinity.

The first result of this section follows from Theorem  \ref{T1} and
Lemma \ref{f02} and describes the situation when all resonances of
$J$ are real and are mostly of the same sign.

\begin{corollary}i) Let $k = 2p$ and $a_p<1.$ Then
$$
\text{\{$J$ has exactly $p$ resonances and all of them are positive\}
$\Leftrightarrow$  \{$\alpha^+_1>2$\}.}
$$

ii) Let $k = 2p$ and $a_p>1, \;\alpha^+_1>2.$

If $\alpha^-_1<-2,$ then $J$ has exactly positive $p-1$ resonances.

If $\alpha^-_1\geq -2$, then $J$ has exactly $p-1$ positive
resonances and 1 negative one.

iii) Let $k = 2p-1.$ Then
$$
\text{\{$J$ has exactly $p-1$
 resonances and all of them are positive and $b_p<0$\}
  $\Leftrightarrow$ \{$\alpha^+_1>2$\}.}
  $$
\end{corollary}
\textbf{Remark.}  Similar arguments can be used when $\a^-_p<-2.$
Then the roots of the Jost function are negative.

\no\textbf{Proof.} i) Let $\alpha^+_1>2.$ Then  there are exactly
$p$ roots $s_1,\dots,s_p$ of $\psi_0$ in $(0,1).$ If $a_p<1,$ then
it follows from Lemma \ref{f01} that there are at least $p$ roots of
$\psi_0$ in $(\frac{1}{s_p}, +\infty).$ Since degree of $\psi_0$ is
equal to $2p$, there are no other roots.

Now, let $J$ have exactly $p$ resonances, all  of them are positive.
Then it follows from Theorem \ref{T1} that $J$ has exactly $p$
eigenvalues. If there is a negative eigenvalue, then it follows from
Lemma \ref{f01} that there is at least one negative resonance. Thus,
all $p$ eigenvalues are positive. Then we can use Theorem \ref{T1}
to obtain that $\alpha^+_j>2, j\in \mathbb{N}_p.$

ii) Let $\alpha^+_1>2$ and $\a^-_1<-2.$ Then it  follows from Lemma
\ref{f01} that there are at least $p-1$ roots of $\psi_0$ in
$(\frac{1}{s_p}, +\infty).$ Since $\psi_0(0)>0$ and
$\psi_0(z)\rightarrow-\infty$ as $z\rightarrow-\infty,$ the last root
has to be in $(-\infty, 0).$ If $\a^-_1<-2,$ then it has to be an
eigenvalue. Thus, there are $1$ negative eigenvalue, $p$ positive
eigenvalues and $p-1$ positive resonances. If $\alpha^-\geq-2,$ then
the last root is a negative resonance and there are $p$ positive
eigenvalues, $p-1$ positive resonances and 1 negative resonance.

iii) If $\alpha^+_1>2,$ then there are exactly  $p$ roots
$s_1,\dots,s_p$ of $\psi_0$ in $(0,1).$  If $b_p>0,$ then it follows
from Lemma \ref{f01} that there are at least $p$ roots of $\psi_0$ in
$(\frac{1}{s_p}, +\infty).$ Since degree of $\psi_0$ is equal to
$2p-1,$ we get a contradiction. Thus, $b_p<0.$ Using the Lemma once
again, we obtain the result.

Now, let $J$ have exactly $p-1$ positive resonances and let $b_p<0.$
Then $J$ has exactly $p$ eigenvalues. If there is a negative
eigenvalue, then it follows from Lemma \ref{f01} that there is at
least one negative resonance. Thus, all eigenvalues are positive, so
$\alpha^+_j>2, j\in \mathbb{N}_p.$  \BBox

\textbf{Proof of Theorems 1.6 and 1.7.} We will  consider the case
$k=2p,$ the proof for the other case is similar. From Lemma
\ref{f01} we know that the product of all roots of $\psi_0(z)$ is
equal to $\frac{1}{1-a^2_p}.$ If $\psi_0(z)$ has no roots in the unit
circle, then the theorem is proved, since the absolute value of all
other roots is greater then 1. Let there be some roots in the real
circle. We already proved that all of them are real and simple.

Let $a_p<1$. Consider the positive roots  $0<s_1<\dots<s_{n_+}<1$ of
$\p_0(z).$ Lemma \ref{f02} shows us that there is at least one root
of $\psi_0(z)$ in every interval $(\frac{1}{s_i},
\frac{1}{s_{i-1}})$ for $i \in[2: n_+].$ In addition, there is a
root in the interval $(\frac{1}{s_1}, +\infty).$ The product of
these roots with all $s_i$ is greater than 1. Doing the same with
negative roots, we obtain that the absolute value of the product of
all roots of $\p_0(z)$ in the real line is greater than 1. From this
we conclude that the absolute value of any resonance outside of the
real line is smaller than $\frac{1}{\sqrt{1-a^2_p}},$ as they always
come in conjugate pairs. Similar arguments work for real resonances
of multiplicity greater than 1. As for the simple real resonances,
the greater of them is the one that is located in the interval
$(\frac{1}{s_1}, +\infty).$ Therefore, the product of all other
resonances is no less than $\frac{1}{s_1\sqrt{1-a^2_p}},$ which
gives us the desired estimate.

If $a_p>1,$ then we can not guarantee that there  is a real root of
$\psi_0(z)$ that is greater than $\frac{1}{s_1},$ similar for
negative roots. Therefore, we can only state that the product of the
resonances outside of the real line is smaller than
$\frac{1}{s_{-1}s_1\sqrt{a^2_p-1}},$ which gives us the estimate for
non-real roots and real roots of multiplicity greater than 1. Note
that there may be no real roots, for example, in $(-1,0).$ Then
$s_{-1}$ is replaced by 1. With simple real roots, there are two
cases. Firstly, the greatest real root can lie in $(1, s_1^{-1}).$
If it is greater than $s_1^{-1},$ then there is another root like
that, as it follows from lemma \ref{f02} that there is an even
amount of resonances in $(s_1^{-1}, +\infty).$ Therefore, the
greatest root has to be smaller than
$\frac{1}{s_{-1}\sqrt{a^2_p-1}}.$ \BBox

\begin{theorem}\lb{corneg}Let $\psi_0$ be the Jost function for some  operator $J$ defined by
(\ref{j}) with a perturbation $q\in \gX_{k}$. Then all its negative
resonances of multiplicity 1 belong to the interval $(-R_-, -1),$
where $R_-$ is defined by:

\no if $k=2p \; (i.e., a_p\neq 1),$ then
\[
R_- = \frac{1}{\beta_-|1 - a_p^2|},\qqq b_- =\ca
\max\{\frac{|s_{-1}|}{a_p^2-1}, s_{1}\}, & {\rm if } \qq a_p<1
\\
|s_1|, & {\rm if }\qq  a_p>1
 \ac,
\]
if $k=2p-1\; (i.e.,  a_p = 1, b_p\neq0),$ then

\[ R_- = \frac{1}{\b_-|b_p|}, \;\;\b_- = \left\{
                                                                 \begin{array}{ll}
                                                                 \max\{\frac{|s_{-1}|}{|b_p|}, 1\},& \hbox{$\text{if } b_p<0$} \\
                                                                  |s_1s_{-1}|, & \hbox{$\text{if } b_p>0$}
                                                                 \end{array}
                                                               \right.
.
\]
\end{theorem}
\no \textbf{Proof.} The proof is similar to the proof of Theorem
1.6. \BBox

\textbf{Proof of Corollary 1.8.} We start with considering the case
$k = 2p, \;a_p<1.$ If $\psi_0(z)$  has exactly $n_+$ positive eigenvalues $s_1<\dots<s_{n_+}$, then there are at least $n_+-1$ roots of
$\psi_0(z)$ in $(s_{n_+}^{-1},s_1^{-1})$, and there can be exactly
$n_+-1$ roots, one in each interval $(s_{j+1}^{-1},s_j^{-1}).$ There
also has to be an odd number of roots in $(s_j^{-1}, +\infty)$,
which gives us at least $2n_+$ roots in $(0, +\infty).$ After
considering the negative half-line we obtain at least $2\gN$ roots
on the real line, which means that we can have at most $p - \gN$
pairs of conjugate complex roots left. If we add real roots instead,
we can obtain the maximal multiplicity $2p - 2\gN + 1$ of the root
in any one of the intervals $(s_{j+1}^{-1},s_j^{-1}).$

Consider the case $a_p>1$ and let both $n_+$ and $n_-$  be greater
than $0$. Then we can have one less root on every half line (it
follows from (\ref{cpinf})), which allows us to have $p - \gN + 1$
pairs of conjugate complex roots or a real root with multiplicity
$2p - 2n + 3.$ If $n_+$ is equal to $0,$ then there is still at
least one positive resonance, since the function $\psi_0(z)$ has to
change its sign between $0$ and $+\infty.$ Therefore, the estimate
is similar to the case when $n_+=1.$ Similar arguments are applied
to $n_-,$ and therefore we obtain the estimate.

Now, consider the case $k = 2p-1.$ Let both $n_+$ be  greater than
1. Similarly, it follows from Lemma \ref{f01} that if $b_p>0,$ then
$\psi_0(z)$ has at least $2n_+$ roots in $(0, +\infty)$ and at least
$2n_- - 1$ roots in $(-\infty, 0),$ which leaves us with no more
than $p - \gN$ pairs of complex roots or with a real root with
multiplicity $2p - 2\gN + 1.$ If $n_+=0,$ then there is a positive
resonance, as in the case ii), so we can replace $n_+$ with 1. The
case $b_p<0$ gives similar result.\BBox

\textbf{Proof of Corollary \ref{cormax}.}
i) Existence  of the matrices
$\wt J$ and $J_o$ follows from Theorem \ref{T4}. Let $k=2p$ for some $p>0.$ Then it follows from
(\ref{f02-1}) that
\[
\Pi_1^{2p}r_n = \frac{1}{1-a^2_p},\;\;\; \frac{\Pi_1^{2p}r_n}{|r|^2}
= \frac{1}{1-(a^o_{p-1})^2}.
\]
Then for $|\ga|\rightarrow+\infty$ we have
\[
\wt a_p^2=1-\frac{|r|^2}{\Pi_1^{2p}r_n}
\frac{1}{|r_{\bu}|^2}=1-\frac{1-(a^o_{p-1})^2}{|\ga|} \rightarrow 1.
\]
Therefore, we have
\[\wt a_p = \Bigl(1-\frac{1-(a^o_{p-1})^2}{|\ga|}\Bigr)^{1/2}
= 1- \frac{1-(a^o_{p-1})^2}{2|\ga|}+o\Bigl(\frac{1}{|\ga|}\Bigr).
\]
Similarly, for $k=2p-1$ we obtain
\[
\wt b_p = \frac{|r|^2}{\Pi_1^{2p-1}r_n}\frac{1}{|r_{\bu}|^2}  =
\frac{b^o_{p-1}}{|\ga|} \rightarrow 0.
\]
Let $\wt\phi(z)$ be the Jost function of the operator $\wt J$  and
let $\phi_o(z)$ be the Jost function for $J_o.$ Then we have
\[\lb{josts}
\wt \phi(z) = \wt C \psi_0(z)\frac{(z-r_{\bu}) (z-\gs_{\bu})}{(z-r)(z-\gs)},\;\;\; \phi_o(z)  =
\frac{C_o\psi_o(z)}{(z-r)(z-\gs)}
\]
for some $\wt C,
C_o\neq 0.$ Substituting (\ref{josts}) into (\ref{fp1}), we obtain
\[\lb{vp1}
\varphi_{p+1}(\wt q, z)(z^2-1) =
\wt C\Bigl( z^{p+2}\psi_0\frac{(\frac{1}{z}-r_\bu)
(\frac{1}{z}-\gs_\bu)}{(\frac{1}{z}-r)
(\frac{1}{z}-\gs)} - \psi_0(z)
\frac{(z-r_\bu)(z-\gs_{\bu})}{z^p(z-r)
(z-\gs)}\Bigr),
\]
\[
\varphi_{p+1}(q_o, z)(z^2-1) = C_o\Bigl( z^{p+2}
\frac{\psi_0(\frac{1}{z})}{(\frac{1}{z}-r)
(\frac{1}{z}-\gs)} -
\frac{\psi_0(z)}{z^p(z-r)(z-\gs)}\Bigr).
\]
Recall that the functions $\varphi_{p+1}(q_o, \l(z))$ and
$\varphi_{p+1}(q_o, \l(z))$ are both polynomials of degree $p.$
After dividing (\ref{vp1}) by $|\ga|$ we obtain that all
coefficients of the polynomial
\[
\frac{\varphi_{p+1}(\wt q, \l(z))}{\wt C |\ga|}
-\frac{\varphi_{p+1}(q_o, \l(z))}{C_o}
\]
tend to zero as $|\ga|$
tends to infinity. Therefore, the mixed eigenvalues $(\wt
\mu_j)_0^{p-1}$ of the matrix $\wt J_p$ tend to those of the matrix
$(J_o)_{p-1}.$ The proof for the Dirichlet eigenvalues $(\wt \t_j)_1^{p-1}$
is similar. Now we can use the continuity of the mapping
$q\rightarrow\mu\star\t$ to obtain that $\wt q\rightarrow q_o$ as
$|\ga|\rightarrow\infty.$

ii) We will prove the existence of the operator $\wt{J},$ the rest of the proof is similar to the proof of the previous case. Let $1< j< n_+.$ If we move the eigenvalue $s_j$ out of the unit circle, we need to have an odd number of resonances in the interval $(s_{j+1}^{-1}, s_{j-1}^{-1})$ to satisfy R4. Since there was an odd number of resonances both in $(s_{j+1}^{-1}, s_{j}^{-1})$ and $(s_{j}^{-1}, s_{j-1}^{-1}),$ then we have an even number in total. If we move one of these resonances out of the interval, this condition is satisfied. As we cannot break other conditions, we need to move both points out of the interval $[s_{-1}^{-1}, s_1^{-1}],$ or to the complex plane, preserving the conjugacy. The proof for the case $j = n_+$ is similar.
\BBox

We can also formulate an analogue of Corollary \ref{cormax}  for
real resonances.

\begin{corollary}\lb{cormax1}
Let the operator $J$ with $q\in \mathfrak{X}_{k}$ be  defined by
(\ref{j}) and let $(r_n)_{1}^{k}$ be the set of roots of its Jost
function.

i) Let $r$ be some resonance lying in the interval $(s_1^{-1}, +\infty)$, where $s_{1}$
is a smallest positive eigenvalues of $J$, or let $r = s_{1}.$ Then $r$ can be moved to any point $\wt r \in \R\setminus[s_{-1}^{-1}, s_1^{-1}]$,
and there exists a matrix $\wt J$ with $\wt q\in\mathfrak{X}_{k}$
corresponding to the obtained set of roots. There also exists a
matrix $J_o$ with $q_0=(a^o, b^o)\in\mathfrak{X}_{k-1},$ which Jost
function has the roots $\{r_n\}_1^{k}\setminus\{r^o\}.$

ii) If $k = 2p\; (a_p\neq1),$ then for $|\wt r|\rightarrow  +\infty$
we have
\[
a_p = 1 - \frac{1-(a^o_{p-1})^2}{2\wt r}+o\Bigl(\frac{1}{\wt r}\Bigr).
\]
If $k = 2p-1\; (a_p = 1, b_p \neq 0),$ then  for $\wt r\rightarrow
+\infty$ we have
\[
b_p = \frac{b_{p-1}}{\wt r} \rightarrow 0.
\]
Moreover, in both cases $\wt q\rightarrow q_o$  as $\wt
r\rightarrow\infty.$
Similar results hold true for negative $r.$
\end{corollary}
\no\textbf{Proof.} The proof is similar to the proof of Corollary
\ref{cormax}. \BBox

\section{Examples}
\setcounter{equation}{0}

{\bf Example 1. }  Consider the case $k = a_1=1, b_1\ne 0$.  Then we
have
\[
\p_0=-b_1z+1,\qqq r_1=1/b_1.
\]
If $|b_1|<1$, then $r_1$ is a resonance. If $|b_1|>1$,  then $r_1$
is an eigenvalue. If $|b_1|=1$, then $r_1$ is a virtual state.

{\bf Example 2.}  Consider the case $k=2, b_1\neq0, a_1\ne 1$.  Then
we have
\[
\p_0={z^{2}c_1-b_1z+1\/a_1} = {c_1\/a_1}(z-r_1)(z-r_2),\qq
r_1r_2={1\/c_1}, \qqq r_2+r_1={b_1\/c_1},
\]
where
\[
r_{1,2}=\frac{b_1\pm \sqrt D}{2c_1},\qqq D=b_1^2-4c_1,\qqq
c_1=1-a_1^2.
\]
There are three cases:

$1_{\bu}$ If $D<0$, then $r_{1,2}$ are complex resonances.

If $r_1\in \C_+\sm \dD_1$, then $r_2=\ol r_1$ and
$c_1={1\/|r_1|^2}<1,\qq b_1={2\Re r_1\/|r_1|^2}.$

$2_{\bu}$ $D=0,$ then $r_1 = r_2 = {b_1\/2c_1}$  is a real resonance
of multiplicity 2.

$3_{\bu}$ If $D>0,$ then $r_{1,2}$ are real.

Moreover, in the last case we obtain:

1) If $r_1,r_2>1$, then $c_1={1\/r_1r_2}>0,
b_1={r_1+r_2\/c_1}>0$.

2) If $r_1,r_2<-1$, then $c_1={1\/r_1r_2}>0,
b_1={r_1+r_2\/c_1}<0$.

3) If $r_1<-1, r_2>1$, then $c_1={1\/r_1r_2}<0,
b_1={r_1+r_2\/c_1}$.

4) If $r_1<-1, r_2\in (0,1)$, then $c_1<0, b_1>0$.

5) If $r_1\in (-1,0), r_2>1$, then $c_1<0, b_1>0$.

6) The cases  $r_1, r_2\in (-1,0)$ or  $r_1, r_2\in (0,1) $  are
impossible, since we have Theorem \ref{T1}.

{\bf Example 3.} Consider a discrete Schrodinger operator  $J_h =
J_0 + hV_p,$ where $h\neq 0,\;$ $J_0$ is a discrete Laplacian on
$\mathbb{N}$ and $V$ is a step  potential given by
\[
(V f)_x= \ca f_x,\;\; x\leq p\\
0, \;\; x>p\ac .\] The matrix $J$ has a perturbation
$q\in\gX_{2p-1}, \; p\in \mathbb{N}. $ Direct calculations show that
for this matrix we  have  $\vp_n(\l)=U_{n-1}(\frac{\l-h}{2})$ and
\[
\textstyle \psi_0(z) = z^pU_{p}(\frac{\l(z)-h}{2})
 - z^{p+1}U_{p-1}(\frac{\l(z)-h}{2}),
 \]
where $U_n(\l)$ are Chebyshev polynomials of the second kind. Then
we obtain
\[
\textstyle
 \mu_0= h-2\cos \Bigl(\frac{\pi }{p+1}\Bigr),\qqq \mu_j =
h-2\text{cos}\Bigl(\pi\frac{j+1}{p+1}\Bigr), \qq j\in \N_{p-1},
\]
\[
\lb{ex31} \textstyle
\alpha^+_j= h-2\cos \Bigl(\frac{2\pi j }
{2p+1}\Bigr),\qqq \alpha^-_j = h-2\text{cos} \Bigl(\frac{-\pi+ 2\pi
j}{2p+1}\Bigr), \qq j\in \N_{p},
\]
and $|h - \m_n|\leq 2, \; n\in [0:p-1].$ The bound states and resonances are located as follows.

1. If $h>2+2\cos \Bigl(\frac{2\pi}{2p+1}\Bigr),$ then
 the matrix $J$ has $p$ eigenvalues $(\l_n)_1^p$, and
\[\lb{ex32}
2<\l_1<\alpha^+_1<\l_{2}<\dots<\l_p<\alpha^+_p.
\]
There are $p$ eigenvalues $(s_n)_1^p, \; s_n = z(\l_n).$
 All $p-1$ resonances $(r_n)_1^{p-1}$ are real
and simple, and $r_i \in (s_{i+1}^{-1}, s_i^{-1}).$  If h tends to
infinity, then it follows from (\ref{ex31}) and (\ref{ex32}) that
\[ \textstyle
s_j = \frac{1+o(1)}{h}, \;\; j\in \mathbb{N}_p, \;\;\; r_j = h+o(h),
\;\; j\in\mathbb{N}_{p-1}.
\]
Note that if we take $(Vf)_x=\epsilon_x f_x,\;\; x\leq p,$ for some
$\epsilon_x\in \R$, then $\frac{1}{h}J_h = \frac{1}{h}J_0 + V_p,$
where $\frac{1}{h}J_0 \rightarrow 0$
 as $h\rightarrow+\infty.$ Therefore, in this case we have
\[
s_j = \frac{1+o(1)}{h\epsilon_j}, \;\; j\in \mathbb{N}_p, \;\;\;
r_j = h\epsilon_j+o(h), \;\; j\in\mathbb{N}_{p-1}.
\]

2. If $h$ decreases to $h = 2+2\cos \Bigl(\frac{2\pi}{2p+1}\Bigr)$,
then $s_p$ tends to 1. After that, $s_p$ turns into a new resonance.
Direct calculations show that it moves to the right, while $r_{p-1}$
moves to the left, until they collide at the point $r^\bu$ at $h =
h_p^\bu$ and become a resonance of multiplicity 2.

3. While  $h_p^\bu>h>2+2\cos \Bigl(\frac{4\pi}{2p+1}\Bigr),$  the
operator $J$ has $p-1$ eigenvalues $(\l_n)_1^{p-1}$, and
\[2<\l_{p-1}<\alpha
_2<
\l_{p-2}<\dots<\l_1<\alpha^+_p.
\]
There are $p-1$ bound states $(s_n)_1^{p-1}, \; s_n = z(\l_n).$
There are $p-3$ resonances $(r_n)_1^{p-3}$ that are still real and
simple. It follows from Theorem \ref{T1} that $r_i \in
(s_{j+1}^{-1}, s_j^{-1}).$ Two other resonances $r_{p-2}$ and
$r_{p-3} = \overline{r_{p-2}}$ are located near $r^\bu.$

4. While $h$ is decreasing further, the largest bound state
$s_{p-1}$ grows until becomes a resonance and then grows some more
before it colliding with $r_{p-3}$ at $h = h^\bu_{p-1},$ generating
two complex resonances. After that, the procedure similar to the one
described in 2-3 repeats again with the bound state $s_{p-2}$ and so
on, until there are no bound states left.

5. If $h$ tends to zero, then all resonances tend to  infinity. In
particular, we have
\[
|r_j| = \frac{1}{2\sqrt[2p-1]{|4h|}}+o(1) \text{ as }
h\rightarrow0, \;\; j\in \mathbb{N}_{2p-1}.
\]
Take
$\frac{\l(z)-h}{2} = \cos\theta.$ Then the points $r^\bu, h_n^\bu$
are defined by the following equations.
\[
r^\bu = \frac{\sin(p+1)\theta}{\sin p\theta},
\]
\[
r^\bu\sin p\theta =  \frac{(r^\bu)^2 - 1}{2}((p+1)\cos(p+1)\theta -
p\cos p\theta).
\]
Note that $\psi_0(z) = 0$ is equivalent to
\[\lb{z}
z = \frac{\sin(p+1)\theta}{\sin p\theta}.\]
If we take the inverse of (\ref{z}) and sum it with (\ref{z}), we obtain
\[
\l(z) = 2\cos\theta + h = \frac{\sin(p+1)
\theta}{\sin p\theta} + \frac{\sin(p)\theta}{\sin (p+1)\theta}.
\]
Therefore,
\[h = \frac{\sin^2(p+1)\theta+ \sin^2(p)\theta - 2\sin(p)\theta \sin(p+1)\theta \cos\theta}{\sin(p+1)\theta \sin(p)\theta} = \frac{\sin^2\theta}{\sin(p+1)\theta \sin(p)\theta}.
\]
If we take $\cos \theta = \eta,$ we obtain
\[\lb{Y}
h = \frac{1 - \eta^2}{\eta - T_{2p+1}(\eta)}.
\]
This gives us $2p-1$ simple roots, $m$ roots in the unit  circle
(the number $m$ depends on the location of the numbers $\alpha _n$,
see Theorem \ref{T1}), $m-1$ roots on the real line outside of it
and $p-m$ pairs of complex roots. If $\eta$ tends to 1, then $\eta -
T_{p+1}(\eta)$ tends to zero, and the fraction in the right part of
(\ref{Y}) tends to some $C\neq 0.$ Since $2p-1$ is an odd number,
the situation when $\eta\rightarrow -1$ is similar. Thus, if $h$
tends to zero, then we have $|\eta|\rightarrow+\infty.$ In
particular, $|\eta| = \frac{1}{2\sqrt[2p-1]{|4h|}}+o(1),$ which
implies
\[|r_n| = \frac{1}{2\sqrt[2p-1]{|4h|}}+o(1) \;\;\text{ as }
\;\; h\rightarrow0.
\]

{\bf Example 4.} Let again $J_h = \vk J_0 + h V$, where $h\neq 0,
\vk>0, \vk\neq 1.$  This matrix has a perturbation $q\in\gX_{2p-1}.$
Then we have $\wt\vp_n(\l)=U_{n-1}(\frac{\l-h}{2\vk}),$ where
$U_n(\l),$ and
\[ \textstyle
\mu_0= h-2\vk\cos \Bigl(\frac{\pi }{p+1}\Bigr),\qqq
 \mu_j = h-2\vk\text{cos}\Bigl(\pi\frac{j+1}{p+1}\Bigr), \qq j\in \N_{p-1},
\]
\[
\textstyle \alpha^+_j= h-2\vk\cos \Bigl(\frac{2\pi j }
{2p+1}\Bigr),\qqq \alpha^-_j = h-2\vk\text{cos} \Bigl(\frac{-\pi+
2\pi j}{2p+1}\Bigr), \qq j\in \N_{p}.
\]
The dynamics of the eigenvalues when $h$ is changed  is similar to
the previous example. The change of $\vk$ leads to the change of
distances between $\alpha^+_j, \m_j$ and, therefore, between the
eigenvalues of the matrix $J$:
\[
|\l_{j+1} - \l_{j}|\rightarrow \infty\;\; \text{ as }\;\;
|\vk|\rightarrow\infty,\;\;\; j\in n_{\bu}\setminus n_+,\]
\[|\l_{j+1} - \l_{j}|\rightarrow 0\;\; \text{ as } \;\;
|\vk|\rightarrow0,\;\;\; j\in n_{\bu}\setminus n_+.
\]

{\bf Example 5.} Consider a matrix $J$ with a perturbation
$q(a,b)\in\gX_{2p},\; p\in \mathbb{N},$ such that $\qqq b_n = 0,
\; n\in\mathbb{N}_p.$ Below we prove that its Jost function
$\psi_0$ satisfies $\psi_0(-z) = \psi_0(z),\; z\in \mathbb{C},$
which means that both eigenvalues and resonances of $J$ come in
pairs, symmetric with respect to zero.

In fact, we prove a stronger statement, namely, that  for such
operators we have
\[\lb{ex5}
\psi_{2x}(-z) = \psi_{2x}(z), \;\; \psi_{2x+1}(-z)
 = -\psi_{2x+1}(z),\;\; \forall \ (x,z)\in \N\ts\mathbb{C}.
 \]
We prove it by induction. We know that for $x>p$ we have $\psi_x(z)
= z^x,$ which gives is the base.  Now, let (\ref{ex5}) be proved for
$x\geq n.$ Then we can use \er{jej} to calculate
\[
\psi_{2n-1}(-z) = \frac{(-z-\frac{1}{z})\psi_{2n}(-z)
 - a_p\psi_{2n+1}(-z)}{a_{p-1}}.
 \]
Since $\psi_{2n}(-z) = \psi_{2n}(z)$ and $\psi_{2n+1}(z)
 = -\psi_{2n+1}(z),$ we obtain that $\psi_{2n-1}(z)
 = -\psi_{2n-1}(z).$ The next step, which is to prove
  that $\psi_{2n-2}(-z) = \psi_{2n-2}(z),$ is similar. Hence,
   (\ref{ex5}) is proved for all $x\in\mathbb{Z}.$

Conversely, let there be some $j$ for which $b_j\neq 0.$ Take  the
maximal $j$ out of these numbers. Then it follows from the first
formula in (\ref{f01-3}) that
$\p_{p-j}(z)={z^{p-j}\/A_{p-j}}\rt(1-zb_j+O(z^{2})\rt)\;\;$ as
$\z\to0.$ Therefore, $\p_{p-j}(-z)\neq\pm\p_{p-j}(z),$ so the
inverse of this statement is also true.

{\bf Example 6. The algortithm for Theorem \ref{pr1}.} The following example demonstrates the algorithm for Theorem \ref{pr1} that was described in the proof of Proposition \ref{TB}. Let the operator $J$ have a perturbation $q\in X^+_p.$ Then it follows from Theorem \ref{pr1} that there is a bijection between $X_p^+$ and $\cE_{2p}(2)$ defined by (\ref{pr1w}). Namely, that we can recover the perturbation $q$ from the sequence of numbers $\omega = (\omega_n)_1^{2p}$ which solve the equation $\frac{2\x(z(\l+0i))}{2p+1} +
\arg z(\l+0i) =  -\frac{\pi n}{2p+1},\; n\in\mathbb{N}_{2p}$. Consider the example of the recovering process in the case $p = 2.$ Let $\omega = (\omega_j)_1^4 = (-1, -\frac12, \frac12, 1).$ Then it follows from the formulas (\ref{TB3}-\ref{TB4}) that we have
$$\alpha_1^- = -1,\;\; \alpha_2^- =\frac12,\;\;\;\;\;\;\alpha_1^+ = -\frac{1}{2},\;\; \alpha_2^+ = 1,
$$
$$a_2 = \Bigl(\frac{(\alpha_1^+ + \alpha_2^+) - (\alpha_1^+ + \alpha_2^+)}{2}\Bigr)^{\frac12} = \frac{\sqrt2}{2},
$$
$$\mathcal{D}^-(\l) = (\l+1)(\l - \frac12), \;\;\;\;\mathcal{D}^+(\l) = (\l+\frac12)(\l-1),
$$
$$\mathcal{D}_{2,1}(\l) = \frac{\mathcal{D}^+(\l) - \mathcal{D^-}(\l)}{2a_2^2} = -\l,\;\;\;\; (\mathcal{D}^+)' = 2\l-\frac12.
$$
Consider the matrix $J^+.$ It follows from (\ref{Tb1}) that
$$(J^+ \delta_2, \delta_2) =-\frac {\alpha_1^+ \mathcal{D}_{2,1}(\alpha^+_1)} {(\mathcal{D}^+)'(\alpha^+_1)} - \frac {\alpha_2^+ \mathcal{D}_{2,1}(\alpha^+_2)}
{(\mathcal{D}^+)'(\alpha^+_2)} = \frac12.$$
On the other hand, direct calculations show that
$$(J^+ \delta_2, \delta_2) = b_2 - a_2^2 = b_2 - \frac12.$$
Therefore, we have $b_2 - \frac12 = \frac12$ and $b_2 = 1.$ Proceeding like this, we obtain that
$$((J^+)^2 \delta_2, \delta_2) =-\frac {(\alpha_1^+)^2 \mathcal{D}_{2,1}(\alpha^+_1)} {(\mathcal{D}^+)'(\alpha^+_1)} - \frac {(\alpha_2^+)^2 \mathcal{D}_{2,1}(\alpha^+_2)}
{(\mathcal{D}^+)'(\alpha^+_2)} = \frac34,$$
while direct calculations show that
$$((J^+)^2 \delta_2, \delta_2) = a_1^2 + (b_2 - a_2)^2 =  a_1^2+\frac14.$$
Since $a_1>0,$ it is determined, and we have $a_1 = \frac{\sqrt2}{2}.$ Finally, we have
$$((J^+)^3 \delta_2, \delta_2) =-\frac {(\alpha_1^+)^3 \mathcal{D}_{2,1}(\alpha^+_1)} {(\mathcal{D}^+)'(\alpha^+_1)} - \frac {(\alpha_2^+)^3 \mathcal{D}_{2,1}(\alpha^+_2)}
{(\mathcal{D}^+)'(\alpha^+_2)} = \frac58,$$
while direct calculations show that
$$((J^+)^3 \delta_2, \delta_2) = b_1a_1^2 +2b_2a_1^2 - 2a_1^2a_2^2 + (b_2 - a_2)^3 = \frac{b_1}{2}+\frac58,
$$
from which we can conclude that $b_1 = 0.$ Therefore, all coefficients of $J$ are found.

\textbf{Remark.} One can see that it takes some time to recover the perturbation $q$ from the numbers $(\omega_j)_1^{2p}$ even in the case $p = 2.$ However, such operations are easily performed using a computer.
\footnotesize\footnotesize

%\no {\bf Acknowledgments.} \footnotesize   Our study was supported
%by the RSF grant No 18-11-00032.

\end{document}